\documentclass[10pt]{article}
\usepackage{latexsym}
\usepackage{epic}
\usepackage{eepic}
\usepackage{color}
\usepackage{pst-node}
\usepackage{amsmath,multirow,amsfonts}
\usepackage{tikz-cd}
\usepackage{tikz}
\usetikzlibrary{matrix}
\usepackage{array}
\usepackage[english]{babel}
\usepackage[latin1]{inputenc}
\usepackage{times}
\usepackage[T1]{fontenc}
\usepackage{amssymb}
\usepackage{graphics}
\usepackage{graphicx}
\usepackage{tikz-cd}
\usepackage{epsfig}
\usepackage{cite}
\usepackage{caption}
\usetikzlibrary{arrows}
\tikzstyle{block}=[draw opacity=2,line width=3cm]
\begin{document}
\newtheorem{def1}{Definition}[section]
\newtheorem{lem}{Lemma}[section]
\newtheorem{exa}{Example}[section]
\newtheorem{thm}{Theorem}[section]
\newtheorem{pro}{Proposition}[section]
\newtheorem{cor}{Corollary}[section]
\newtheorem{rem}{Remark}[section]
\newtheorem{exam}{Example}[section]
\title{On unification of categories associated with $F$-transforms and fuzzy pretopological spaces as Qua category}
\author{Abha Tripathi\thanks{tripathiabha29@gmail.com}\\
Department of Mathematics, School of Basic Sciences\\
Galgotias University\\
Greater Noida 203201 , U.P., India \\ and\\
S.P. Tiwari\thanks{sptiwari.1maths@gmail.com}\\
Department of Mathematics \& Computing\\ Indian Institute of Technology (ISM)\\
Dhanbad-826004, India}
\date{}
\maketitle
\begin{abstract}
In this contribution, our motive is to unify the categories associated with $F$-transforms and fuzzy pretopological spaces as a new category ${\bf Qua}$, whose object classes are success measurements of answers and morphisms are pairs of success measurements of transformations. Specifically, the categories of spaces with $L$-valued fuzzy partitions, $L$-valued fuzzy lower transformation systems, $L$-valued fuzzy pretopological spaces, and \v{C}ech $L$-valued fuzzy interior spaces with the morphisms as pairs of $L$-valued fuzzy relations between the underlying sets of corresponding objects have intriguing relationships with the category ${\bf Qua}$.
\end{abstract}
\section{Introduction} 
The concept of fuzzy transform ($F$-transform) was firstly introduced by Perfilieva \cite{per}, a notion that piqued the curiosity of many researchers. It has now been greatly expanded upon, and a new chapter in the theory of semi-linear spaces has been opened. The  fundamental idea of the $F$-transform is to factorize (or fuzzify) the precise values of independent variables by using a proximity relationship, and to average the precise values of dependent variables to an approximation value (cf., \cite{{per},{irin1}}), from fuzzy sets to parametrized fuzzy sets \cite{st} and from the single variable to the two (or more variables) (cf., \cite{{ma},{mar}, {mar1}, {step}}). Recently, several studies have begun to look into $F$-transforms based on an arbitrary $L$-fuzzy partition of an arbitrary universe (cf., \cite{kh1,mockor1,jir, mocko, anan,spt1}), where $L$ is a complete residuated lattice. Among these researches, the concept of a general transformation operator defined by a monadic relation was introduced in \cite{mockor1}, the links between $F$-transforms and {semimodule homomorphisms} were investigated in \cite{jir}, while the connections between $F$-transforms and similarity {relations} were discussed in \cite{mocko}. Further, a fascinating relationship was also discovered in \cite{anan} among $F$-transforms, $L$-fuzzy topologies/co-topologies, and $L$-fuzzy approximation operators (all of which are ideas employed in the study of an operator-oriented perspective of rough set theory), while the relationships of $L^M$-valued fuzzy rough sets and $ML$-graded topologies/co-topologies with $L^M$-valued $F$-transforms was discussed in \cite{abhac,trip}. Also, a connection between fuzzy pretopological spaces and spaces with $L$-fuzzy partition was investigated in \cite{spt1}. In which it has been shown that weaker closure and interior operators, called after \v{C}ech, may also be expressed by using $F$-transforms, implying that lattice-based $F$-transforms can be utilized in parallel with closure and interior operators as their canonical representation. Furthermore, classes of $F$-transforms taking into account three well-known classes of implicators, namely {$R-,S-, QL-$ implicators} were discussed in \cite{tri}. Additionally, a new notion of $F$-transforms based on overlap and grouping maps, residual, and co-residual implicator over complete lattice from constructive and axiomatic approaches have been investigated in \cite{abhaog}. Several studies in the subject of $F$-transforms applications have been conducted, e.g., trend-cycle estimation \cite{holc}, {data compression \cite{hut}}, numerical solution of partial differential equations \cite{kh,abha}, scheduling \cite{li}, time series \cite{vil}, data analysis \cite{no}, denoising \cite{Ir}, face recognition \cite{roh}, neural network { approaches \cite{ste} and trading} \cite{to}.\\\\
The notion of the category theory initiated by Eilenberg and Mac Lane \cite{eil} is well-known, which was further evolved by many scholars in \cite{fre,gro,law,law1,pi}. In the framework of $F$-transforms, the category ${\bf SpaceFP}$ was introduced by Mo\v{c}ko\v{r} \cite{mock,jiri,jir,mocko,mockor} as a generalization of categories of sets with fuzzy partitions determined by lattice valued fuzzy equivalences (or, equivalently, similarity relations). In \cite{mock,jiri,jiri1,mockor}, the category ${\bf SpaceFP}$ of spaces with fuzzy partitions and some properties have been investigated. Also, in \cite{jiri}, it has been shown that ${\bf SpaceFP}$ is a topological category, and is isomorphic to the category of upper (lower) transformation systems that satisfy simplified axioms. Moreover, in \cite{mock}, it has been shown that there exists a functorial relationship among a subcategory of ${\bf SpaceFP}$ and categories of Kuratowski closure and interior operators, and a category of fuzzy preorders, respectively. {In \cite{jiri1}, for generalized powerset theories in categories, new conceptions of relational, closure, or partition powerset theories in these categories were studied.} Further, the functorial relationships among the categories of sets with \v{C}ech closure and interior operators, $L$-fuzzy pretopologies and $L$-fuzzy co-pretopologies, reflexive $L$-fuzzy relations, and $L$-fuzzy partitions have been studied in \cite{mockor}. Also, in \cite{qi}, the concepts of $L$-fuzzy pretopological spaces and $L$-fuzzy approximation spaces based on the reflexive $L$-fuzzy relations from a categorical viewpoint have been investigated. In a different direction, many categories of $L$-fuzzy sets have been demonstrated to be special instances of the category ${\bf Qua}$ in various ways in \cite{sinha}. Specifically, the category of $L$-fuzzy sets defined in \cite{gog}, category ${\bf RelL}$ defined in \cite{bar1}, the category ${\bf M_L( Set)}$ defined in \cite{pai}, the category ${\bf PV}$ defined in \cite{blas} and the category {\bf GameK} defined in \cite{laf} are particular cases of the category ${\bf Qua}$. 
Motivated from the above, the unification of the categories of spaces with $L$-valued fuzzy partitions, $L$-valued fuzzy transformation systems, $L$-valued fuzzy pretopological spaces, and \v{C}ech $L$-valued fuzzy interior spaces is the objective of this work. Specifically, in this paper, we present a categorical approach in the framework of spaces with $L$-valued fuzzy partitions, $L$-valued fuzzy lower transformation systems, $L$-valued fuzzy pretopological spaces and \v{C}ech $L$-valued fuzzy interior spaces with morphisms as pair of $L$-valued fuzzy relations and establish an isomorphisms between the categories of spaces with $L$-valued fuzzy partitions and $L$-valued fuzzy pretopological spaces, $L$-valued fuzzy lower transformation systems and \v{C}ech $L$-valued fuzzy interior spaces, respectively. We particularize the above-introduced categories in numerous ways by utilizing the notion of the category ${\bf Qua}$. Furthermore, we introduce and study adjoint functors among the studied categories. The main findings are summarized below:
\begin{itemize}
\item we introduce the categories ${\bf LSpaceFP}$ (category of spaces with $L$-valued fuzzy partitions), ${\bf LFtrans}^\downarrow$ (category of $L$-valued fuzzy lower transformation systems), ${\bf LFPrTop}$ (category of $L$-valued fuzzy pretopological spaces) and ${\bf LFCInt}$ (category of \v{C}ech $L$-valued fuzzy interior spaces) with morphisms as pairs of $L$-valued fuzzy relations;
\item we demonstrate that the categories ${\bf LSpaceFP}$ and ${\bf LFPrTop}$ are isomorphic to the categories ${\bf LFtrans}^\downarrow$ and ${\bf LFCInt}$, respectively;
\item we establish a functorial relationship among the categories ${\bf LSpaceFP}$, ${\bf LFtrans}^\downarrow$, ${\bf LFPrTop}$  and ${\bf LFCInt}$; 
\item we show that there exist two pairs of adjoint functors between the categories ${\bf LSpaceFP}$, ${\bf LFPrTop}$ and the categories ${\bf LFtrans}^\downarrow$, ${\bf LFCInt}$, respectively; and
\item we show the relationships of the categories associated with $F$-transforms and $L$-valued fuzzy pretopological spaces with the category {\bf Qua}.
\end{itemize} 
\section{Preliminaries}
Herein, we collect basic ideas related to the category theory and complete lattices. There are two subsections in this section.
\subsection{Category theory}
\begin{figure}
\[\begin{tikzcd}[row sep=12ex, column sep=12ex] 
   C\arrow{r}{\Phi_C} \arrow[swap]{dr}{f} & G(F(C))\arrow{d}{G(g)} &F(C)\arrow{d}{g}\\
     & G(D)&D
  \end{tikzcd}\]
 \caption{Diagram for Definition \ref{adj0}.}
\label{fig:0}
  \end{figure}
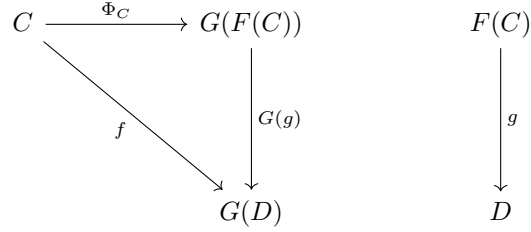

{We collect some well-known ideas relating to categories, subcategories, functors, which are to be used in the main text. For details on category theory, we refer \cite{ad,ar,bar,fre,lan}.} Throughout, for a category ${\bf C}$, $|{\bf C}|$ denotes the class of objects of ${\bf C}$, while its morphisms are written as ${\bf C}$-morphisms. We begin with the following.
\begin{def1}
{A category ${\bf D}$ is said to be a} {\bf subcategory} of a category ${\bf C}$ if
\begin{itemize}
\item[(i)] every ${\bf D}$-object is a ${\bf C}$-object ,
\item[(ii)] for all ${\bf D}$-objects $D$ and $D'$, ${\bf D}(D, D')\subseteq{\bf C}(D, D')$, and
\item[(iii)] {category ${\bf D}$ has the same composition and identity morphisms as category ${\bf C}$}. 
\end{itemize}
\end{def1}
Assume that ${\bf C}$ and ${\bf D}$ are two categories for the rest of this subsection. Now, we have the following.
\begin{def1}
{A function $F : {\bf C} \rightarrow {\bf D}$ is said to be a functor if it maps every ${\bf C}$-object $C$ to a ${\bf D}$-object $F(C)$ and every ${\bf C}$-morphism $\phi : C \rightarrow C'$ to a ${\bf D}$-morphism $F(\phi) : F(C) \rightarrow F(C')$ such that
\begin{itemize}
\item[(i)] for each ${\bf C}$-morphisms $\phi : C \rightarrow C'$ and $\phi': C' \rightarrow C''$, $F(\phi' \circ \phi ) = F(\phi') \circ F(\phi)$, and
\item[(ii)] for every ${\bf C}$-objects, $F(id_C) = id_{F(C)}$.
\end{itemize}}
\end{def1}
\begin{def1} Two categories ${\bf C}$ and ${\bf D}$ are isomorphic, if there exist two functors  $F:{\bf C}\rightarrow {\bf D}$ and $ G:{\bf D}\rightarrow {\bf C}$ such that $G\circ F=id_{\bf C}$ and $F\circ G=id_{\bf D}$
\end{def1}
\begin{def1}\label{adj0}
{For two functors $F:{\bf C}\rightarrow{\bf D},G:{\bf D}\rightarrow {\bf C}$ be functors, if there exists a natural transformation $\Phi:id_{\bf C}\rightarrow{\bf G\circ F}$ such that for each ${\bf C}$-object $C$, ${\bf D}$-object $D$ and ${\bf C}$-morphism $\phi:C\rightarrow G(D)$, there exists a unique  ${\bf D}$-morphism $\psi:F(C)\rightarrow D$ such that the diagram in Figure \ref{fig:0} commutes. Then $F$ is {\bf left adjoint} to $G$ and $G$ is {\bf right adjoint} to $F$.} 
\end{def1}
  \subsection{Complete residuated and co-residuated lattice}
{Herein, we collect some basic concepts} based on complete residuated and co-residuated lattices from  \cite{belo,gog,kl,kli,qia}. Besides, we recall the concept of the category ${\bf Qua}$ whose object-class are success measurements of answers and morphisms are pairs of $L$-valued fuzzy relations from \cite{sinha}. We begin with the definition of a complete residuated lattice.
\begin{def1}
{A {\bf complete residuated lattice} is an algebra $(L,\wedge,\vee,\ast,0,1)$ such that}
  \begin{itemize}  
  \item[(i)] $(L,\wedge,\vee,0,1)$ is a complete lattice with the least element $0$ and the greatest element $1$;
  \item[(ii)] $(L,\ast,1)$ is a commutative monoid;
  \item[(iii)] $a\ast \bigvee\limits_{i\in J}b_i=\bigvee\limits_{i\in J}a\ast b_i,\,\forall\, \{b_i:i\in J\}\subseteq L$.
  \end{itemize}
  \end{def1}
  \begin{def1}
A {\bf complete co-residuated lattice} is an algebra $(L,\wedge,\vee,\#,0,1)$ suchh that
  \begin{itemize}  
  \item[(i)] $(L,\wedge,\vee,0,1)$ is a complete lattice with the least element 0 and the greatest element 1;
  \item[(ii)] $(L,\#,0)$ is a commutative monoid;
  \item[(iii)] $a\# \bigwedge\limits_{i\in J}b_i=\bigwedge\limits_{i\in J}a\# b_i,\,\forall\, \{b_i:i\in J\}\subseteq L$.
  \end{itemize}
  \end{def1}
  \begin{def1}
A {\bf negator} on $L$ is a decreasing function $\neg:L\rightarrow L$ such that $\neg 0 =1,\neg1=0$. Also, $\neg$ is said to be an {\bf involutive negator} if $\neg\neg a=a,\,\forall\,a\in L$.
  \end{def1} 
Throughout this work, we consider that $L\equiv(L,\wedge,\vee,\ast,\#,0,1)$ is a complete residuated and co-residuated lattice and $\neg$ is an involutive negator, in which $\ast$ and $\#$ associated by absorptin law, i.e., $(a\#b)\ast c\leq a\#(b\ast c),\,\forall\,a,b,c\in L$ and dual with respect $\neg$, i.e., $\neg a\ast \neg b=\neg(a\#b),\,\forall\,a,b\in L$. For an involutive negator $\neg$, $\bigwedge\limits_{i\in J}(\neg a_i)=\neg(\bigvee\limits_{i\in J} a_i)$ and $\bigvee\limits_{i\in J}(\neg a_i)=\neg(\bigwedge\limits_{i\in J} a_i),\,\forall\, \{a_i:i\in J\}\subseteq L.$\\\\
For $X \in |{\bf SET}|$, an {\bf $L$-valued fuzzy} set in $X$ is a function $f:X\rightarrow L$ and the collection of all $L$-valued fuzzy sets in $X$ is denoted by $L^X$. Also, for all $x\in X,a\in L,\textbf{{a}}\in L^X,\textbf{{a}}(x)=a$ denotes constant $L$-valued fuzzy set. 
\begin{rem}
The {\bf core} of an $L$-valued fuzzy set $f$ is defined as: 
\begin{eqnarray*}
core(f)=\lbrace x\in X:f(x)=1\rbrace.
\end{eqnarray*}
If { $core(f)\neq \emptyset$}, then $f$ is called {a \textbf {normal} $L$-valued fuzzy} set. For $A\subseteq X$, the {\bf characteristic function} of $A$ is a function $1_A:X\rightarrow \{0,1\}$ such that
\begin{eqnarray*} 
1_A(x)=\begin{cases}
1 &\text{ if } x\in A,\\
0&\text{ otherwise}.
\end{cases} 
\end{eqnarray*}
\end{rem}
\begin{rem}
 For $f,g\in L^X$, $f\leq g$ if $f(x)\leq g(x),\,\forall\,x\in X$. Also, for a given system $\lbrace f_i:i\in J\rbrace\subseteq L^X$,
 $$(\bigwedge\limits_{i\in J}f_i )(x)=\bigwedge\limits_{i\in J} f_i(x) ~\text{and} ~(\bigvee\limits_{i\in J}f_i)(x)=\bigvee\limits_{i\in J} f_i(x),\forall x\in X.$$
 \end{rem}
 Now, let $\alpha : X \rightarrow Y$ be a $|{\bf SET}|$-morphism. Then according to Zadeh's extension principle, $\alpha$ can be extended to the forward and backward operators ${\alpha}^\rightarrow : L^{X} \rightarrow L^{Y}$ and ${\alpha}^\leftarrow : L^{Y}\rightarrow L^{X}$ such that for all $y\in {Y},f \in L^{X}, g \in L^{Y}$
$${\alpha}^\rightarrow(f)(y)=\bigvee\limits_{x\in X,\alpha(x)=y}f(x),\,\,{\alpha}^\leftarrow(g)=g\circ \alpha,\,\text{respectively}.$$
{Also, for $X,Y,Z\in |{\bf SET}|$, let $R:X\times Y\rightarrow L$ and $R':Y\times Z\rightarrow L$ be $L$-valued fuzzy relations. Then the multiplication of $R$ and $R'$ over complete lattice $(L,\wedge,\vee,\ast,0,1)$ is an $L$-valued fuzzy relation $R'\odot R:X\times Z\rightarrow L$ such that
\[(R'\odot R)(x,z)=\bigvee\limits_{y_1\in Y_1}R(x,y)\ast R'(y,z),\,\forall\,x\in X,z\in Z,\,\text{and}\]
the multiplication of $R$ and $R'$ over complete co-residuated lattice $(L,\wedge,\vee,\#,0,1)$ is an $L$-valued fuzzy relation $R'\oplus R:X\times Z\rightarrow L$ such that
\[(R'\oplus R)(x,z)=\bigwedge\limits_{y_1\in Y_1}R(x,y)\# R'(y,z),\,\forall\,x\in X,z\in Z.\]}{
The following are towards the induced operations of intersection $\wedge$, union $\vee$, multiplications $\ast$,$\#$ and negator $\neg$ on $L^X$.} For all $x\in X,f,g\in L^X$
\begin{itemize}
    \item[(i)] $f=g\Leftrightarrow f(x)=g(x),\,f\leq g\Leftrightarrow f(x)\leq g(x)$,
    \item[(ii)] $(f\wedge g)(x)=f(x)\wedge g(x),\,(f\vee g)(x)=f(x)\vee g(x)$,
    \item[(iii)] $(f\ast g)(x)=f(x)\ast g(x),\,(f\# g)(x)=f(x)\# g(x)$, and
    \item[(iv)] $(\neg f)(x)=\neg f(x)$.
\end{itemize}
Next, we recollect the concepts of the category ${\bf Qua}$ from \cite{sinha}.
\begin{def1} {The category {\bf Qua}} is defined as follows:\\\\
{\bf Object:} A system $\mathcal{M}=(M^{*},M_{*},\mathbf{M})$, where $M^{*}$ is a set whose members are referred to as ``{\bf  type $M$ questions}'', $M_{*}$ is a set whose members are referred to as ``type $M$ {\bf answers}'' and $\mathbf{M}:M^{*}\times M_{*}\rightarrow L$ is an $L$-valued fuzzy relation such that for each $m^{*}\in M^{*}$ and $m_{*}\in M_{*},\,\mathbf{M}(m^{*},m_{*})\in L$ is referred to as the ``{\bf success-measure}" of answer $m_{*}$ in response to the question $m^{*}$.\\\\
{\bf Morphism:} : For objets $\mathcal{M}=(M^{*},M_{*},\mathbf{M})$ and $\mathcal{M}'=(M'^{*},M'_*,\mathbf{M}')$, the morphism from $\mathcal{M}$ to $\mathcal{M}'$ is a pair $(g_*,g^*)$ such that
\begin{itemize}
    \item[(i)] $g_*:M_{*}\times M'_*\rightarrow L$ and $g^*:M'^{*}\times M^{*}\rightarrow L$ are $L$-valued fuzzy relations
    \item[(ii)] $M_*(m_{*},m'_*)\ast \mathbf{M}(m^{*},m_{*})\leq g^*(m'^{*},m^{*})\#\mathbf{M'}(m'^{*},m'_*),\,\forall\,m_{*}\in M_{*},m'_*\in M'_*,m^{*}\in M^{*},m'^{*}\in M'^{*}$.
\end{itemize}
{\bf Identity:} For object $\mathcal{M}=(M^{*},M_{*},\mathbf{M})$, the identity morphism on $\mathcal{M}$ is a pair of $L$-valued fuzzy relations ${id_{M_{*}}}:M_{*}\times M_{*}\rightarrow L$ and $id_{M^{*}}:M^{*}\times M^{*}\rightarrow L$ such that $id_{M_{*}}=1_{M_{*}}$ and ${id_{M^{*}}}=0_{M^{*}}.$ \\\\
{\bf Composition:}  For objects $\mathcal{M}=(M^{*},M_{*},\mathbf{M})$, $\mathcal{M}'=(M'^{*},M'_*,\mathbf{M}')$ and $\mathcal{M}''=(M''^{*},M''_*,\mathbf{M}'')$ and morphisms $(g_*,g^*):\mathcal{M}\rightarrow\mathcal{M}'$ and $(f'_*,f'^{*}):\mathcal{M}'\rightarrow\mathcal{M}''
$, the {\bf composition} of $(g_*,g^*)$ and $(g'_*,g'^{*})$ is a pair of $L$-valued fuzzy relations $(g'_*,g'^{*})\bullet (g_*,g^*):\mathcal{M}\rightarrow \mathcal{M}''$ such that $(g'_*,g'^{*})\bullet (g_*,g^*)=(g'_*\odot g_*,g^*\oplus g'^{*})$.
\end{def1}
\begin{rem}
\begin{itemize}
\item[(i)] {The following is the motivation for this particular morphism condition.\\\\
Let $\mathcal{M},\mathcal{M}'\in {\bf Qua}$ and $(g_*,g^*):\mathcal{M}\rightarrow \mathcal{M}'$ be ${\bf Qua}$-morphisms. Furthermore, let $m_{*}$ be a type $M$ answer in response to the question $m^{*}$ and success of this response is measured by $\mathbf{M}(m^{*}, m_{*})$. At the same time, $g_*(m_{*},m'_*)$ measures the success of transformation of an answer $m_{*}$ into the  type $M'$ answer $m'_*$. The whole success of this operation is determined by $g_*(m_{*}, m'_*)\ast \mathbf{M}(m^{*},m_{*})$, which is the result of the cooperation of $g_*(m_{*}, m'_*)$ and $ \mathbf{M}(m^{*},m_{*})$. On the other hand, if $m'^{*}$ is a type $M'$ question, the answer $m'_*$ will react to it and the success of this response is determined by $\mathbf{M}'(m'^{*}, m'_*)$. At the same time, $g^{*}(m'^{*}, m^{*})$ measures the success of transformation of a question $m'^{*}$ into the  type $M$ question $m^{*}$. The whole success of this operation is determined by $g^{*}(m'^{*}, m^{*})\#\mathbf{M}'(m'^{*}, m'_*)$, which is the result of the competition between $g^{*}(m'^{*}, m^{*})$ and $\mathbf{M}'(m'^{*}, m'_*)$. Thus the defining condition of $L$-valued fuzzy homomorphism is that the result of cooperation of $g_*(m_{*}, m'_*)$ and $ \mathbf{M}(m^{*},m_{*})$ should not grater than the
result of the competition between $M^{*}(m'^{*}, m^{*})$ and $\mathbf{M}'(m'^{*}, m'_*)$.
\item[(ii)] The class of object of the category ${\bf Qua}$ is clearly success measurements of responses and morphisms are pairs of success measurements of transformations.}
\end{itemize}
\end{rem}
\section{Categories of spaces with $L$-valued fuzzy partitions and $L$-valued fuzzy lower transformation systems}
Herein, we discuss the ideas of the category of spaces with $L$-valued fuzzy partitions with morphisms as pairs of $L$-valued fuzzy relations, as well as the category of $L$-valued fuzzy lower transformation systems. Furthermore, we establish a connection of the above-introduce categories with ${\bf Qua}$. Now, we recall the following from \cite{anan}.
\begin{def1}\label{FP}
A collection $\mathcal{P}$ of normal $L$-valued fuzzy sets $\lbrace A_{j}:j\in J\rbrace$ is called an {\bf $L$-valued fuzzy partition} of $X$ if the corresponding collection of ordinary sets $\lbrace core(A_{j}):j\in J\rbrace$ is partition of $X$. The pair $(X,\mathcal{P})$ is called a {\bf space with $L$-valued fuzzy partition}.
\end{def1}
For any $L$-valued fuzzy partition $\mathcal{P} = \lbrace A_{j}:j\in J\rbrace$, it is possible to associate it with the following onto index map $\xi:X\rightarrow J$ such that
$$\xi(x)=j\Leftrightarrow x\in core(A_{j}),\,\forall\,x\in X,j\in J.$$
{In the following, the concept of the category spaces with $L$-valued fuzzy partitions is introduced and their morphisms are defined as pairs of $L$-valued fuzzy relations between the underlying sets of corresponding objects. Now, we start with the following.}
\begin{def1}
Let ${\pi}_1=(X_1,J_1,\mathcal{P}_1)$ and $\pi_2=(X_2,J_2,\mathcal{P}_2)$ be spaces with $L$-valued fuzzy partitions, where $\mathcal{P}_1 = \lbrace A_{j_1}:j_1\in J_1\rbrace$ and $\mathcal{P}_2 = \lbrace A_{j_2}:j_2\in J_2\rbrace$. Then an {\bf $L$-valued fuzzy FP-map} from ${\pi}_1$ to $\pi_2$ is a pair $(\phi,\mu)$ such that
\begin{itemize}
\item[(i)] $\phi:X_1\times X_2\rightarrow L$ and $\mu:J_2\times J_1\rightarrow L$ are $L$-valued fuzzy relations, and
\item[(ii)]  $\phi(x_1,x_2)\ast A_{j_1}(x_1)\leq \mu(j_2,j_1)\#A_{j_2}(x_2),\,\forall\,x_1\in X_1,x_2\in X_2,j_1\in J_1,j_2\in J_2$.
\end{itemize}
\end{def1}
Let ${\pi}_1=(X_1,J_1,\mathcal{P}_1),\pi_2=(X_2,J_2,\mathcal{P}_2) ,\pi_3=(X_3,J_3,\mathcal{P}_3)$ be spaces with $L$-valued fuzzy partitions and $(\phi,\mu):\pi_1\rightarrow \pi_2,(\phi',\mu'):\pi_2\rightarrow \pi_3 $ be $L$-valued fuzzy FP-maps. Then the {\bf composition} $(\phi',\mu')\bullet  (\phi,\mu):\pi_1\rightarrow \pi_3 $ of $(\phi,\mu)$ and $(\phi',\mu') $ is a pair of $L$-valued fuzzy relations such that $ (\phi',\mu')\bullet (\phi,\mu)=(\phi'\odot \phi,\mu\oplus\mu')$.
\begin{pro}\label{1} Spaces with $L$-valued fuzzy partitions with their $L$-valued fuzzy FP-maps form a category.
\end{pro}
\textbf{Proof:} It is enough to demonstrate that the composition of two $L$-valued fuzzy FP-maps is also an $L$-valued fuzzy FP-map. To do this, let $\pi_1=(X_1,J_1,\mathcal{P}_1),\pi_2=(X_2,J_2,\mathcal{P}_2),\pi_3=(X_3,J_3,\mathcal{P}_3) $ be spaces with $L$-valued fuzzy partitions and $(\phi,\mu):\pi_1\rightarrow \pi_2,(\phi',\mu'):\pi_2\rightarrow \pi_3 $ be $L$-valued fuzzy FP-maps. Then for all $x_1\in X_1,x_2\in X_2,j_1\in J_1$
\begin{eqnarray*}
(\phi'\odot\phi)(x_1,x_3)\ast A_{j_1}(x_1)&=&\bigvee\limits_{x_2\in X_2}(\phi(x_1,x_2)\ast\phi'(x_2,x_3))\ast A_{j_1}(x_1)\\
&=&\bigvee\limits_{x_2\in X_2}\phi'(x_2,x_3)\ast (\phi(x_1,x_2)\ast A_{j_1}(x_1))\\
&\leq&\bigvee\limits_{x_2\in X_2}\phi'(x_2,x_3)\ast (\mu(j_2,j_1)\# A_{j_2}(x_2))\\
&\leq&\bigvee\limits_{x_2\in X_2}\mu(j_2,j_1)\#(\phi'(x_2,x_3)\ast A_{j_2}(x_2))\\
&\leq&\mu(j_2,j_1)\#(\mu'(j_3,j_2)\# A_{j_3}(x_3))\\
&\leq&\bigwedge\limits_{j_2\in J_2}(\mu(j_2,j_1)\#\mu'(j_3,j_2))\# A_{j_3}(x_3)\\
&=&(\mu\oplus \mu')(j_3,j_1)\# A_{j_3}(x_3).
\end{eqnarray*}
Thus $(\phi'\odot \phi,\mu\oplus\mu'):\pi_1\rightarrow\pi_2$ is an $L$-valued fuzzy FP-map.\\\\
We shall denote by ${\bf LSpaceFP}$, the category of spaces with $L$-valued fuzzy partitions with $L$-valued fuzzy FP-maps as morphisms.\\\\
In the following, we show that there exists a subcategory of the category ${\bf LSpaceFP}$, which is isomorphic to the category of spaces with $L$-valued fuzzy partitions (say, {\bf SpaceFP}) introduced in \cite{mock}, whose morphisms are pairs of functions between the underlying sets of corresponding objects. For which, we assume that ${\bf SLSpaceFP}$ is a category, whose objects are same as the objects of ${\bf LSpaceFP}$ and morphisms are ${\bf LSpaceFP}$-morphisms $(\phi,\mu):\pi_1\rightarrow \pi_2$ such that
\begin{itemize}
\item[$\bullet$] $\phi:X_1\times X_2\rightarrow L$ is an $L$-valued fuzzy relation and for each $x_1\in X_1$ there exists a unique $x_2\in X_2, \phi(x_1,x_2)=1$, and 0, otherwise,
\item[$\bullet$] $\mu:J_2\times J_1\rightarrow L$ is an $L$-valued fuzzy relation and for each $j_1\in J_1,$ there exists a unique $j_2\in J_2, \mu(j_2,j_1)=0$ and 1, otherwise.
\end{itemize}
{It is obvious that ${\bf SLSpaceFP}$ is a subcategory of} ${\bf LSpaceFP}$.
\begin{pro}\label{2} Let $F_1:{\bf SpaceFP}\rightarrow {\bf SLSpaceFP}$ be a function such that for every $\pi_1\in |{\bf SpaceFP}|,F_1(\pi_1)=\pi_1$ and for every ${\bf SpaceFP}$-morphism $(\alpha,\omega):\pi_1\rightarrow \pi_2,F_1(\alpha,\omega):\pi_1\rightarrow \pi_2$ be a function such that $F_1(\alpha,\omega)=(\phi_{\alpha},\mu_{\omega})$, where $\phi_{\alpha}:X_1\times X_2\rightarrow L,\mu_{\omega}:J_2\times J_1\rightarrow L$ are $L$-valued fuzzy relations such that for each $x_1\in X_1,x_2\in X_2,j_1\in J_1,j'\in J'$,
\begin{eqnarray*} 
\phi_{\alpha}(x_1,x_2)=\begin{cases}
1 &\text{ if } \alpha(x_1)=x_2\\
0&\text{ otherwise,} 
\end{cases}
 \mu_{\omega}(j_2,j_1)=\begin{cases}
0 &\text{ if } \omega(j_1)=j_2\\
1&\text{ otherwise,}
\end{cases}
\end{eqnarray*}
respectively. Then $F_1$ is a functor.
\end{pro}
\textbf{Proof:} {It is easy to understand the proof.}
\begin{pro}\label{3} Let $F'_1:{\bf SLSpaceFP}\rightarrow {\bf SpaceFP}$ be a function such that for every $\pi_1\in |{\bf SLSpaceFP}|,F'_1(\pi_1)=\pi_2$ and for every ${\bf SLSpaceFP}$-morphism $(\phi,\mu):\pi_1\rightarrow \pi_2,F'_1(\phi,\mu):\pi_1\rightarrow \pi_2$ be a function such that $F_1'(\phi,\mu)=({\alpha}_{\phi},\omega_{\mu})$, where $\alpha_{\phi}:X_1\rightarrow X_2$ and $\omega_{\mu}:J_1\rightarrow J_2$ are function such that for all $x_1\in X_1,x_2\in X_2,j_1\in J_1,j_2\in J_2,\,\alpha_{\phi}(x_1)=x_2$ iff $\phi(x_1,x_2)=1$ and $\omega_{\mu}(j_1)=j_2$ iff $\mu(j_2,j_1)=0$. Then $F_1'$ is a functor. 
\end{pro}
\textbf{Proof:} {It is easy to understand the proof.}
\begin{pro}\label{4} The categories ${\bf SpaceFP}$ and ${\bf SLSpaceFP}$ are isomorphic.
\end{pro}
\textbf{Proof:} Let $F_1:{\bf SpaceFP}\rightarrow {\bf SLSpaceFP}$ and $F'_1:{\bf SLSpaceFP}\rightarrow {\bf SpaceFP}$ be functors. Then $F'_1\circ F_1:{\bf SpaceFP}\rightarrow {\bf SpaceFP}$ is a functor such that for all $\pi_1\in|{\bf SpaceFP}|$ and ${\bf SpaceFP}$-morphism $(\alpha,\omega):\pi_1\rightarrow\pi_2$, 
\[(F'_1\circ F_1)(\pi_1)=F'_1( F_1(\pi_1))=F'_1(\pi_1)=\pi_1.\]
Thus $(F'_1\circ F_1)(\pi_1)=\pi_1$ and 
\[(F'_1\circ F_1)(\alpha,\omega)=F'_1( F_1(\alpha,\omega))=F'_1(\phi_{\alpha},\mu_{\omega})=(\alpha_{\phi_\alpha},\omega_{\mu_\omega})=(\alpha,\omega).\]
Thus $(F'_1\circ F_1)(\alpha,\omega)=(\alpha,\omega)$ and therefore $F'_1\circ F_1=id_{{\bf SpaceFP}}$. Similarly, we can show that $F_1\circ F'_1:{{\bf SLSpaceFP}}\rightarrow {{\bf SLSpaceFP}}$ is a functor such that for all $\pi_1\in|{{\bf SLSpaceFP}}|$ and ${{\bf SLSpaceFP}}$-morphism $(\phi,\mu):\pi_1\rightarrow\pi_2$, 
$(F_1\circ F'_1)(\pi_1)=\pi_1,(F_1\circ F'_1)(\phi,\mu)=(\phi,\mu)$. Thus $F_1\circ F'_1=id_{{\bf SLSpaceFP}}.$ Therefore categories ${\bf SpaceFP}$ and ${\bf SLSpaceFP}$ are isomorphic.\\\\
In the following, we unify the categories ${\bf LSpaceFP}$ and ${\bf SpaceFP}$ as ${\bf Qua}$ category.
\begin{pro}\label{5} The category ${\bf LSpaceFP}$ is a subcategory of the category ${\bf Qua}$.
\end{pro}
\textbf{Proof:} Let $\pi_1=(X_1,J_1,\mathcal{P}_1)\in |{\bf LSpaceFP}|$ and $(\phi,\mu):\pi_1\rightarrow \pi_2$ be an ${\bf LSpaceFP}$-morphism. Then $\pi_1=(X_1,J_1,\mathcal{P}_1)$ can be seen as an object $ \mathcal{M}_{\pi_1}=(M^{*},M_{*},\mathbf{M}_{\mathcal{P}_1})\in|{\bf Qua}|$, where $M^{*}=J_1,M{_*}=X_1$ and $\mathbf{M}_{\mathcal{P}_1}:J_1\times X_1\rightarrow L$ is an $L$-valued fuzzy relation such that $\mathbf{M}_{\mathcal{P}_1}(j_1,x_1)=A_{j_1}(x_1),\,\forall\,j_1\in J_1,x_1\in X_1$. We can also construct a pair of functions $(g_{*},g^{*}):\mathcal{M}_{\pi_1}\rightarrow \mathcal{M}_{\pi_1}$ such that $ g_{*}:X_1\times X_2\rightarrow L,g^{*}:J_2\times J_1\rightarrow L$ are $L$-valued fuzzy relations and $g_{*}(x_1,x_2)=\phi(x_1,x_2),g^{*}(j_2,j_1)=\mu(j_2,j_1),\,\forall\,x_1\in X_1,x_2\in X_2,j_1\in J_1,j_2\in J_2$. Now, let $x_1\in X_1,x_2\in X_2,j_1\in J_1$. Then
\begin{eqnarray*}
g_{*}(x_1,x_2)\ast\mathbf{M}_{\mathcal{P}_1}(j_1,x_1)&=&\phi(x_1,x_2)\ast A_{j_1}(x_1)\\
&\leq&\mu(j_2,j_1)\#A_{j_2}(x_2)\\
&=&g^{*}(j_2,j_1)\# \mathbf{M}_{\mathcal{P}_2}(j_2,x_2).
\end{eqnarray*}
Thus $(g_{*},g^{*}):\mathcal{M}_{\pi_1}\rightarrow \mathcal{M}_{\pi_2}$ is a ${\bf Qua}$-morphism. Therefore the category ${\bf LSpaceFP}$ is a subcategory of the category ${\bf Qua}$.
\begin{pro}\label{6} There is a subcategory of the category ${\bf Qua}$, which is isomorphic to the category ${\bf SpaceFP}$.
\end{pro}
\textbf{Proof:} {Propositions \ref{4} and \ref{5} lead to the proof.}\\\\
Next, we introduce the concept of the lower $F$-transform. We begin with the following.
\begin{def1}\label{lt}
Let $\mathcal{P}=\lbrace A_{j}:j\in J\rbrace$ be an $L$-valued fuzzy partition. Then the {\bf direct $F^\downarrow$-transform} of $f\in L^X$ is a is a collection of lattice elements $\{F^\downarrow_j[f]:j\in J\}$ and the $j$-th component of direct $F^\downarrow$-transform is given by
$$F_j^\downarrow[f]=\bigwedge\limits_{x\in X}(\neg A_{j}(x)\# f(x)),\, \forall\,j\in J,f\in L^X.$$
\end{def1}
For $f \in L^X$ and for an $L$-valued fuzzy partition $\mathcal{P}$, it can be seen that the lower $F$-transform $(F^\downarrow)$ induces a function $F^\downarrow_{\mathcal{P}}:L^X\rightarrow L^J$ such that
$F^\downarrow_{\mathcal{P}}[f](j) = F^\downarrow_j[f]$. 
\begin{pro}\label{3.1}
Let $\mathcal{P}=\lbrace A_{j}:j\in J\rbrace$ be an $L$-valued fuzzy partition. Then 
\item[(i)] $F_{\mathcal{P}}^\downarrow[\textbf{a}]=\textbf{a},\,\forall\,a\in L$,
\item[(ii)] $F_{\mathcal{P}}^\downarrow[f](j)\leq f(x_j),\,\forall\,x_j\in core(A_j)$,
\item[(iii)] $F_{\mathcal{P}}^\downarrow[\textbf{a}\# f]=\textbf{a}\# F_{\mathcal{P}}^\downarrow[f],\,\forall\,\textbf{a},f\in L^X$, and
\item[(iv)]  $F_{\mathcal{P}}^\downarrow[\bigwedge\limits_{i\in J}f_i]=\bigwedge\limits_{i\in J}F_{\mathcal{P}}^\downarrow[f_i],\,\forall\,\{f_i:i\in J\}\subseteq L^X$.
\end{pro}
\textbf{Proof:} Definition \ref{lt} leads to this proof.\\\\
Now, we introduce the concept of an $L$-valued fuzzy lower transformation system.
\begin{def1}
Let $X$ be a nonempty set. Then the system $\mathcal{H}=(X,Y,v,H)$, where 
\begin{itemize}
\item[1.] $Y$ is a set,
\item[2.] $v:X\rightarrow Y$ is a surjective function,
\item[3.] $H:L^X\rightarrow L^Y$ is such that
\begin{itemize}
\item[(i)] $H(\bigwedge\limits_{i\in J}f_i)=\bigwedge\limits_{i\in J}H(f_i),\,\forall\,\{f_i:i\in J\}\subseteq L^X,$
\item[(ii)] $H(\textbf{a}\# f)=\textbf{a}\# H(f),\,\forall\,\textbf{a},f\in L^X$, and
\item[(iii)] $\neg H(\neg 1_{\{x\}})(y)=1$ iff $v(x)=y,\,\forall\,x\in X,y_1\in Y_1$,
\end{itemize}
\end{itemize} 
is called an {\bf $L$-valued fuzzy lower transformation system} on $X$.
\end{def1}
\begin{def1}
Let $\mathcal{H}_1=(X_1,Y_1,v_1,H_1)$ and $\mathcal{H}_2=(X_2,Y_2,v_2,H_2)$ be $L$-valued fuzzy lower transformation systems on $X_1$ and $X_2$, respectively. Then the {\bf homomorphism} from $\mathcal{H}_1$ to $\mathcal{H}_2$ is a pair $(\beta,\epsilon)$ such that
\begin{itemize}
\item[(i)] $\beta:X_1\rightarrow  X_2$ and $\epsilon:Y_1\rightarrow  Y_2$ are functions, and
\item[(ii)] ${\epsilon}^\leftarrow( H_2(f))\leq H_1( {\beta}^\leftarrow(f)),\,\forall\,f\in L^{X_2}$.
\end{itemize}
\end{def1}
Let $\mathcal{H}_1=(X_1,Y_1,v_1,H_1),\mathcal{H}_2=(X_2,Y_2,v_2,H_2),\mathcal{H}_3=(X_3,Y_3,v_3,H_3)$ be $L$-valued fuzzy lower transformation systems on $X_1,X_2,X_3$, respectively and $(\beta,\epsilon):\mathcal{H}_1\rightarrow \mathcal{H}_2,(\beta',\epsilon'):\mathcal{H}_2\rightarrow \mathcal{H}_3$ be the homomorphisms. Then the {\bf composition} $(\beta',\epsilon')\circ  (\beta,\epsilon):\mathcal{H}_1\rightarrow \mathcal{H}_3$ of $(\beta,\epsilon)$ and $(\beta',\epsilon') $ is a pair of functions such that $ (\beta',\epsilon')\circ  (\beta,\epsilon)=(\beta'\circ \beta,\epsilon'\circ\epsilon)$.
\begin{pro} $L$-valued fuzzy lower transformation systems with their homomorphisms form a category.
\end{pro}
\textbf{Proof:} {It is easy to understand the proof.}\\\\
We shall denote by ${\bf Ftrans}^{\downarrow}$, the category of $L$-valued fuzzy lower transformation systems with homomorphisms as morphisms.\\\\
In the following, we use the notion of power set functions defined by an $L$-valued fuzzy relation, which was first defined in \cite{gog}.
If $\psi:X_1\times X_2\rightarrow L$ be an $L$-valued fuzzy relation, then the power set function $\overleftarrow{\psi}:L^{X_2}\rightarrow L^{X_1}$ is given as
$$\overleftarrow{\psi}(f)(x_1)=\bigwedge\limits_{x_2\in X_2}(\neg\psi(x_1,x_2)\# f(x_2)),\,\forall\,x_1\in X_1,f\in L^{X_2}.$$
Now, let $\psi:X_1\times X_2\rightarrow L, \psi':X_2\times X_3\rightarrow L$ and $\psi'\odot\psi:X_1\times X_3\rightarrow L$ be $L$-valued fuzzy relations and corresponding to these $L$-valued fuzzy relations we define power set functions $\overleftarrow{\psi}:L^{X_2}\rightarrow L^{X_1},\overleftarrow{\psi}':L^{X_3}\rightarrow L^{X_2}$ and $\overleftarrow{\psi'\odot\psi}:L^{X_3}\rightarrow L^{X_1}$ by 
\[\overleftarrow{\psi}(f)(x_1)=\bigwedge\limits_{x_2\in X_2}\neg\psi(x_1,x_2)\#f(x_2),\,\forall\,x_1\in X_1,f\in L^{X_2},\]
\[\overleftarrow{\psi}'(g)(x_2)=\bigwedge\limits_{x_3\in X_3}\neg\psi'(x_2,x_3)\#g(x_3),\,\forall\,x_2\in X_2,g\in L^{X_3},\, and\]
\[(\overleftarrow{\psi'\odot\psi})(g)(x_1)=\bigwedge\limits_{x_3\in X_3}\neg({\psi'\odot\psi})(x_1,x_3)\#g(x_3),\,\forall\,x_1\in X_1,g\in L^{X_3}.\]
Then $\overleftarrow{\psi'\odot\psi}=\overleftarrow{\psi}\circ\overleftarrow{\psi'}$.\\\\
To show this, let $x_1\in X_1,g\in L^{X_3}$. Then
\begin{eqnarray*}
(\overleftarrow{\psi'\odot\psi})(g)(x_1)&=&\bigwedge\limits_{x_3\in X_3}\neg({\psi'\odot\psi})(x_1,x_3)\#g(x_3)\\
&=&\bigwedge\limits_{x_3\in X_3}\neg(\bigvee\limits_{x_2\in X_2}(\psi(x_1,x_2)\ast\psi'(x_2,x_3))\#g(x_3)\\
&=&\bigwedge\limits_{x_3\in X_3}\bigwedge\limits_{x_2\in X_2}\neg(\psi(x_1,x_2)\ast\psi'(x_2,x_3))\#g(x_3)\\
&=&\bigwedge\limits_{x_3\in X_3}\bigwedge\limits_{x_2\in X_2}(\neg\psi(x_1,x_2)\#\neg\psi'(x_2,x_3))\#g(x_3)\\
&=&\bigwedge\limits_{x_2\in X_2}\neg\psi(x_1,x_2)\#\bigwedge\limits_{x_3\in X_3}\neg\psi'(x_2,x_3)\#g(x_3)\\
&=&\bigwedge\limits_{x_2\in X_2}\neg\psi(x_1,x_2)\#\overleftarrow{\psi}'(g)(x_2)\\
&=&\overleftarrow{\psi}(\overleftarrow{\psi}'(g))(x_1)\\
&=&(\overleftarrow{\psi}\circ \overleftarrow{\psi}')(g)(x_1).
\end{eqnarray*}
Thus $\overleftarrow{\psi'\odot\psi}=\overleftarrow{\psi}\circ \overleftarrow{\psi}'$.\\\\
{In the following, the concept of the category $L$-valued fuzzy lower transformation systems is introduced and their morphisms are defined as pairs of $L$-valued fuzzy relations between the underlying sets of corresponding objects.}
\begin{def1}
Let $\mathcal{H}_1=(X_1,Y_1,v_1,H_1)$ and $\mathcal{H}_2=(X_2,Y_2,v_2,H_2)$ be $L$-valued fuzzy lower transformation systems on $X_1$ and $X_2$, respectively. Then an {\bf $L$-valued fuzzy homomorphism} from $\mathcal{H}_1$ to $\mathcal{H}_2$ is a pair $(\psi,\nu)$ such that
\begin{itemize}
\item[(i)] $\psi:X_1\times X_2\rightarrow L$ and $\nu:Y_2\times Y_1\rightarrow L$ are $L$-valued fuzzy relations, and
\item[(ii)] $\neg\nu(y_2,y_1)\ast H_2(f)(y_2)\leq  H_1(\overleftarrow{\psi}(f))(y_1),\,\forall\,y_1\in Y_1,y_2\in Y_1,f\in L^{X_2}$.
\end{itemize}
\end{def1}
Let $\mathcal{H}_1=(X_1,Y_1,v_1,H_1),\mathcal{H}_2=(X_2,Y_2,v_2,H_2),\mathcal{H}_3=(X_3,Y_3,v_3,H_3)$ be $L$-valued fuzzy lower transformation systems on $X_1,X_2,X_3$, respectively and $(\psi,\nu):\mathcal{H}_1\rightarrow \mathcal{H}_2,(\psi',\nu'):\mathcal{H}_2\rightarrow \mathcal{H}_3 $ be their $L$-valued fuzzy homomorphisms. Then the {\bf composition} $(\psi',\nu')\bullet  (\psi,\nu):\mathcal{H}_1\rightarrow \mathcal{H}_3 $ of $ (\psi,\nu)$ and $(\psi',\nu') $ is a pair of $L$-valued fuzzy relations such that $(\psi',\nu')\bullet  (\psi,\nu)=(\psi'\odot \psi,\nu\oplus\nu')$.
\begin{pro}\label{8} $L$-valued fuzzy lower transformation systems with their $L$-valued fuzzy homomorphisms form a category.
\end{pro}
\textbf{Proof:} It is enough to demonstrate that the composition of two $L$-valued fuzzy homomorphisms is also an $L$-valued fuzzy homomorphism. To do this, let $\mathcal{H}_1=(X_1,Y_1,v_1,H_1),\linebreak\mathcal{H}_2=(X_2,Y_2,v_2,H_2),\mathcal{H}_3=(X_3,Y_3,v_3,H_3)$ be $L$-valued fuzzy lower transformation systems and $(\psi,\nu):\mathcal{H}_1\rightarrow \mathcal{H}_2, (\psi',\nu'):\mathcal{H}_2\rightarrow \mathcal{H}_3$ be $L$-valued fuzzy homomorphisms. Then for all $y_1\in Y_1,f\in L^{X_3}$
\begin{eqnarray*}
H_1((\overleftarrow{\psi'\odot \psi})(f))(y_1)&=&H_1((\overleftarrow{\psi}\circ \overleftarrow{\psi'})(f))(y_1)\\
&=&H_1(\overleftarrow{\psi}( \overleftarrow{\psi'}(f)))(y_1)\\
&\geq&\neg\nu(y_2,y_1)\ast H_2( \overleftarrow{\psi'}(f))(y_2)\\
&\geq&\neg\nu(y_2,y_1)\ast \neg\nu'(y_3,y_2)\ast H{_3}( f)(y_3)\\
&=&\neg(\nu(y_2,y_1)\#\nu'(y_3,y_2))\ast H{_3}( f)(y_3)\\
&\geq&\bigvee\limits_{y_2\in Y_2}\neg(\nu(y_2,y_1)\#\nu'(y_3,y_2))\ast H{_3}( f)(y_3)\\
&=&\neg(\bigwedge\limits_{y_2\in Y_2}\nu(y_2,y_1)\#\nu'(y_3,y_2))\ast H{_3}( f)(y_3)\\
&=&\neg(\nu\oplus\nu')(y_3,y_1)\ast H{_3}( f)(y_3).
\end{eqnarray*}
Thus $(\psi'\odot \psi,\nu\oplus\nu'):\mathcal{H}_1\rightarrow\mathcal{H}_3$ is an $L$-valued fuzzy homomorphism.\\\\
We shall denote by ${\bf LFtrans}^{\downarrow}$, the category of $L$-valued fuzzy lower transformation systems with $L$-valued fuzzy homomorphisms as morphisms.\\\\
In the following, we show that there exists a subcategory of the category ${\bf LFtrans}^{\downarrow}$, which is isomorphic to the category of $L$-valued fuzzy lower transformation systems (say, ${\bf Ftrans}^{\downarrow}$), whose morphisms are pairs of functions between the underlying sets of corresponding objects. For which, we assume that ${\bf SLFtrans}^{\downarrow}$ is a category, whose objects are same as the objects of ${\bf LFtrans}^{\downarrow}$ and morphisms are ${\bf LFtrans}^{\downarrow}$-morphisms $(\psi,\nu):\mathcal{H}_1\rightarrow \mathcal{H}_2$ such that 
\begin{itemize}
\item[$\bullet$] $\psi:X_1\times X_2\rightarrow L$ is an $L$-valued fuzzy relation and for each $x_1\in X_1$, there exists a unique $x_2\in X_2, \psi(x_1,x_2)=1$ and 0, otherwise,
\item[$\bullet$] $\nu:Y_2\times Y_1\rightarrow L$ is an $L$-valued fuzzy relation and for each $y_1\in Y_1$, there exists a unique $y_2\in Y_2,\,\mu(y_2,y_1)=0$ and 1, otherwise.
\end{itemize}
{It is obvious that ${\bf SLFtrans}^{\downarrow}$ is a subcategory} of ${\bf LFtrans}^{\downarrow}$.
\begin{pro}\label{9} Let $F_2:{{\bf Ftrans}^{\downarrow}}\rightarrow {{\bf SLFtrans}^{\downarrow}}$ be a function such that for every $\mathcal{H}_1\in |{\bf Ftrans}^{\downarrow}|,F_2(\mathcal{H}_1)=\mathcal{H}_1$ and for every ${\bf Ftrans}^{\downarrow}$-morphism $(\beta,\epsilon):\mathcal{H}_1\rightarrow \mathcal{H}_2,F_2(\beta,\epsilon):\mathcal{H}_1\rightarrow \mathcal{H}_2$ be a function such that $F_2(\beta,\epsilon)=(\psi_{\beta},\nu_{\epsilon})$, where $\psi_{\beta}:X_1\times X_2\rightarrow L,\nu_{\epsilon}:Y_2\times Y_1\rightarrow L$ are $L$-valued fuzzy relations such that for each $x_1\in X_1,x_2\in X_2,y_1\in Y_1,y_2\in Y_2$,
\begin{eqnarray*} 
\psi_{\beta}(x_1,x_2)=\begin{cases}
1 &\text{ if } \beta(x_1)=x_2\\
0&\text{ otherwise,} 
\end{cases}
 \nu_{\epsilon}(y_2,y_1)=\begin{cases}
0 &\text{ if } \epsilon(y_1)=y_2\\
1&\text{ otherwise,}
\end{cases}
\end{eqnarray*}
respectively. Then $F_2$ is a functor.
\end{pro}
\textbf{Proof:} {It is easy to understand the proof.}
\begin{pro}\label{10} Let $F'_2:{{\bf SLFtrans}^{\downarrow}}\rightarrow {{\bf Ftrans}^{\downarrow}}$ be a function such that for every $\mathcal{H}_1\in |{\bf SLFtrans}^{\downarrow}|,F'_2(\mathcal{H}_1)=\mathcal{H}_1$ and for every ${\bf SLFtrans}^{\downarrow}$-morphism $(\psi,\nu):\mathcal{H}_1\rightarrow \mathcal{H}_2,F'_2(\psi,\nu):\mathcal{H}_1\rightarrow \mathcal{H}_2$ be a function such that $F_2'(\psi,\nu)=({\beta}_{\psi},{\epsilon}_\nu)$, where $\beta_{\psi}:X_1\rightarrow X_2$ and ${\epsilon_\nu}:Y_1\rightarrow  Y_2$ are functions such that for each $x_1\in X_1,x_2\in X_2,y_1\in Y_1,y_2\in Y_2,\,\,\beta_{\psi}(x_1)=x_2$ iff $\psi(x_1,x_2)=1$ and $\epsilon_{\nu}(y_1)=y_2$ iff $\nu(y_2,y_1)=0$. Then $F_2'$ is a functor. 
\end{pro}
\textbf{Proof:} {It is easy to understand the proof.}
\begin{pro} \label{0}
The categories ${\bf Ftrans}^\downarrow$ and ${\bf SLFtrans}^\downarrow$ are isomorphic. 
\end{pro}
\textbf{Proof:} Let $F_2:{{\bf Ftrans}^{\downarrow}}\rightarrow {{\bf SLFtrans}^{\downarrow}}$ and $F'_2:{{\bf SLFtrans}^{\downarrow}}\rightarrow {{\bf Ftrans}^{\downarrow}}$ be functors. Then $F'_2\circ F_2:{{\bf Ftrans}^{\downarrow}}\rightarrow {{\bf Ftrans}^{\downarrow}}$ is a functor such that for all $\mathcal{H}_1\in|{{\bf Ftrans}^{\downarrow}}|$ and ${{\bf Ftrans}^{\downarrow}}$-morphism $(\beta,\epsilon):\mathcal{H}_1\rightarrow\mathcal{H}_2$, 
\begin{eqnarray*}
(F'_2\circ F_2)(\mathcal{H}_1)&=&F'_2( F_2(\mathcal{H}_1))\\
&=&F'_2(\mathcal{H}_1)\\
&=&\mathcal{H}_1.
\end{eqnarray*}
Thus $(F'_2\circ F_2)(\mathcal{H}_1)=\mathcal{H}_1$ and 
\begin{eqnarray*}
(F'_2\circ F_2)(\beta,\epsilon)&=&F'_2( F_2(\beta,\epsilon))\\
&=&F'_2(\psi_\beta,\nu_\epsilon)\\
&=&(\beta_{\psi_\beta},\epsilon_{\nu_\epsilon})\\&=&(\beta,\epsilon).
\end{eqnarray*}
Thus $(F'_2\circ F_2)(\beta,\epsilon)=(\beta,\epsilon)$ and therefore $F'_2\circ F_2=id_{{\bf Ftrans}^{\downarrow}}$. Similarly, we can show that $F_2\circ F'_2:{{\bf SLFtrans}^{\downarrow}}\rightarrow {{\bf SLFtrans}^{\downarrow}}$ is a functor such that for all $\mathcal{H}_1\in|{{\bf SLFtrans}^{\downarrow}}|$ and ${{\bf SLFtrans}^{\downarrow}}$-morphism $(\psi,\nu):\mathcal{H}_1\rightarrow\mathcal{H}_2$, $(F_2\circ F'_2)(\mathcal{H}_1)=\mathcal{H}_1,(F_2\circ F'_2)(\psi,\nu)=(\psi,\nu)$.
Thus $F_2\circ F'_2=id_{{\bf SLFtrans}^{\downarrow}}.$
Therefore the categories ${\bf Ftrans}^\downarrow$ and ${\bf SLFtrans}^\downarrow$ are isomorphic.\\\\
In the following, we show that the category ${\bf LSpaceFP}$ is isomorphic to the category ${\bf LFtrans}^\downarrow$. For which, let $\pi=(X,J,\mathcal{P})$ be a space with $L$-valued fuzzy partition and $\mathcal{H}_1=(X,Y,v,H)$ be an $L$-valued fuzzy lower transformation system. Then we can introduce an $L$-valued fuzzy lower transformation system and space with $L$-valued fuzzy partition with respect to $\pi$ and $\mathcal{H}$, respectively, defined by $\mathcal{H}_\pi$ and $\pi_{\mathcal{H}}$, as follows:
 \begin{itemize}
     \item[(i)] $\mathcal{H}_\pi=(X,J,\xi,H_{\mathcal{P}})$, where $\xi:X\rightarrow J$ is a surective function and $H_{\mathcal{P}}:L^X\rightarrow L^J$ such that $
     H_{\mathcal{P}}(f)(j)=\bigwedge\limits_{x_1\in X_1}\neg A_{j}(x)\#f(x)=F^\downarrow_{\mathcal{P}}[f](j),\,\forall\,j\in J,f\in L^X,$ and
     \item[(ii)] $\pi_{\mathcal{H}}=(X,Y,\mathcal{P}_{H})$, where $\mathcal{P}_{H}=\{A^{H}_{y}:y\in Y\}$ and $v:X\rightarrow Y$ is a surjective index function such that for each $x\in X,y\in Y,\,v(x)=y\Leftrightarrow x\in core(A^{H}_y)$, $A^{H}_y(x)=\neg H(\neg 1_{\{x\}})(y)$ and $H=F^\downarrow_{\mathcal{P}_H}$.
 \end{itemize}
\begin{pro}\label{12}
Each ${\bf LSpaceFP}$-morphism is also an ${\bf LFtrans}^\downarrow$-morphism.
\end{pro}
\textbf{Proof:} Let $\pi_1=(X_1,J_1,\mathcal{P}_1),\pi_2=(X_2,J_2,\mathcal{P}_2)\in |{\bf LSpaceFP}|$ and $(\phi,\mu):\pi_1\rightarrow\pi_1$ be an ${\bf LSpaceFP}$-morphism. Then for all $j_1\in J_1,f\in L^{X_2}$
\begin{eqnarray*}
H_{\mathcal{P}_1}(\overleftarrow{\phi}(f))(j_1)&=&F^\downarrow_{\mathcal{P}_1}[\overleftarrow{\phi}(f)](j_1)\\
&=&\bigwedge\limits_{x_1\in X_1}\neg A_{j_1}(x_1)\#\overleftarrow{\phi}(f)(x_1).\\
&=&\bigwedge\limits_{x_1\in X_1}\neg A_{j_1}(x_1)\#\bigwedge\limits_{x_2\in X_2}\neg\phi(x_1,x_2)\#f(x_2)\\
&=&\bigwedge\limits_{x_1\in X_1}\bigwedge\limits_{x_2\in X_2}\neg (A_{j_1}(x_1)\ast\phi(x_1,x_2))\#f(x_2)\\
&\geq&\bigwedge\limits_{x_2\in X_2}\neg (\mu(j_2,j_1)\#A_{j_2}(x_2))\#f(x_2)\\
&=&\bigwedge\limits_{x_2\in X_2}(\neg \mu(j_2,j_1)\ast \neg A_{j_2}(x_2))\#f(x_2)\\
&\geq&\bigwedge\limits_{x_2\in X_2}\neg \mu(j_2,j_1)\ast (\neg A_{j_2}(x_2)\#f(x_2))\\
&\geq&\neg \mu(j_2,j_1)\ast \bigwedge\limits_{x_2\in X_2} (\neg A_{j_2}(x_2)\#f(x_2))\\
&=&\neg \mu(j_2,j_1)\ast F^{\downarrow}_{\mathcal{P}_2}[f](j_2)\\
&=&\neg \mu(j_2,j_1)\ast H_{\mathcal{P}_2}(f)(j_2).
\end{eqnarray*}
Thus $(\phi,\phi'):\mathcal{H}_{\pi_1}\rightarrow \mathcal{H}_{\pi_2}$ is an ${\bf LFtrans}^\downarrow$-morphism.
\begin{pro}\label{13}
Each ${\bf LFtrans}^\downarrow$-morphism is also an ${\bf LSpaceFP}$-morphism.
\end{pro}
\textbf{Proof:} Let $\mathcal{H}_1=(X_1,Y_1,v_1,H_1),\mathcal{H}_2=(X_2,Y_2,v_2,H_2)\in |{\bf LFtrans}^\downarrow|$ and $(\psi,\nu):\mathcal{H}_1\rightarrow\mathcal{H}_2$ be an ${\bf LFtrans}^\downarrow$-morphism. Then for all $y_1\in Y_1,x_2\in X_2$
\begin{eqnarray*}
 H_1(\overleftarrow{\psi}(\neg1_{\{x_2\}}))(y_1)
&=& H_1(\bigwedge\limits_{x_1\in X_1}\overleftarrow{\psi}(\neg1_{\{x_2\}})(x_1)\#\neg1_{\{x_1\}})(y_1)\\
&=& \bigwedge\limits_{x_1\in X_1}\overleftarrow{\psi}(\neg1_{\{x_2\}})(x_1)\#H_1(\neg1_{\{x_1\}})(y_1)\\
&=&\bigwedge\limits_{x_1\in X_1}(\bigwedge\limits_{x_2'\in X_2}\neg\psi(x_1,x_2')\#\neg1_{\{x_2
\}}(x_2'))\#H_1(\neg1_{\{x_1\}})(y_1)\\
&=&\bigwedge\limits_{x_1\in X_1}\neg\psi(x_1,x_2)\#H_1(\neg1_{\{x_1\}})(y_1)\\
&\leq&\neg\psi(x_1,x_2)\#H(\neg1_{\{x_1\}})(y_1)\\
&=&\neg(\psi(x_1,x_2)\ast\neg H_1(\neg1_{\{x_1\}})(y_1))\\
&=&\neg(\psi(x_1,x_2)\ast A^{{H_1}}_{y_1}(x_1)).
\end{eqnarray*}
On the other hand, for all $y_1\in Y_1,y_2\in Y_2,x_2\in X_2$ 
\begin{eqnarray*}
\nu(y_2,y_1)\#A^{{H}{_2}}_{y_2}(x_2)&=&\nu(y_2,y_1)\#\neg H_1(\neg1_{\{x_2\}})(y_2)\\
&=&\neg(\neg\nu(y_2,y_1)\ast H_2(\neg1_{\{x_2\}})(y_2))\\
&\geq&\neg( H_1(\overleftarrow{\psi}(\neg1_{\{x_2\}}))(y_1))\\
&\geq&\neg\neg(\psi(x-1,x_2)\ast A^{{H_1}}_{y_1}(x_1))\\
&=&\psi(x_1,x_2)\ast A^{{H_1}}_{y_1}(x_1).
\end{eqnarray*}
Thus $(\psi,\nu):\pi_{\mathcal{H}_1}\rightarrow\pi_{\mathcal{H}_2}$ is an ${\bf LSpaceFP}$-morphism.
\begin{pro}\label{14}
Let $F_3:{\bf LSpaceFP}\rightarrow {\bf LFtrans}^\downarrow$ be a function such that for every $\pi_1\in|{\bf LSpaceFP}|,F_3(\pi_1)=\mathcal{H}_{\pi_1}$ and for every ${\bf LSpaceFP}$-morphism $(\phi,\mu):\pi_1\rightarrow\pi_2,F_3(\phi,\mu):\mathcal{H}_{\pi_1}\rightarrow\mathcal{H}_{\pi_2}$ be a function such that $F_3(\phi,\mu)=(\phi,\mu)$. Then $F_3$ is a functor.
\end{pro}
\textbf{Proof:} (i) Let $(\phi,\mu):\pi_1\rightarrow\pi_2$ and $(\phi',\mu'):\pi_2\rightarrow\pi_3$ be ${\bf LSpaceFP}$-morphisms. Then 
\begin{eqnarray*}
F_3((\phi{'},\mu')\bullet(\phi,\mu))&=&(\phi{'},\mu')\bullet(\phi,\mu)\\
&=&F_3(\phi{'},\mu')\bullet F_3(\phi,\mu).
\end{eqnarray*}
Thus $F_3((\phi{'},\mu')\bullet(\phi,\mu))=F_3(\phi{'},\mu')\bullet F_3(\phi,\mu)$.\\\\
(ii) Let $\pi_1\in |{\bf LSpaceFP}|$. Then $F_3(id_{\pi_1})=F_3(1_{X_1},0_{J_1})=(1_{X_1},0_{J_1})$. Thus $F_3(id_{\pi_1})=id_{F_3(\pi_1)}$. Therefore $F_3$ is a functor.
\begin{pro}\label{15}
Let $F_3':{\bf LFtrans}^\downarrow\rightarrow {\bf LSpaceFP}$ be a function such that for every $\mathcal{H}_1\in|{\bf LFtrans}^\downarrow|,F_3(\mathcal{H}_1)=\pi_{\mathcal{H}_1}$ and for every ${\bf LFtrans}^\downarrow$-morphism $(\psi,\nu):{\mathcal{H}_1}\rightarrow{{\mathcal{H}_2}},F_3'(\psi,\nu):\pi_{\mathcal{H}_1}\rightarrow\pi_{{\mathcal{H}_2}}$ be a function such that $F_3'(\psi,\nu)=(\psi,\nu)$. Then $F_3'$ is a functor.
\end{pro}
\textbf{Proof:} (i) Let $(\psi,\nu):\mathcal{H}_1\rightarrow\mathcal{H}_2$ and $(\psi',\nu'):\mathcal{H}_2\rightarrow\mathcal{H}_3$ be ${\bf LFtrans}^{\downarrow}$-morphisms. Then 
\begin{eqnarray*}
F'_3((\psi{'},\nu')\bullet(\psi,\nu))&=&(\psi{'},\nu')\bullet(\psi,\nu)\\
&=&F'_3(\psi{'},\nu')\bullet F'_3(\psi,\nu).
\end{eqnarray*}
Thus $F'_3((\psi{'},\nu')\bullet(\psi,\nu))=F'_3(\psi{'},\nu')\bullet F'_3(\psi,\nu)$.\\\\
(ii) Let $\mathcal{H}_1\in |{\bf LFtran}^{\downarrow}|$. Then $F'_3(id_{\mathcal{H}_1})=F'_3(1_{X_1},0_{Y_1})=(1_{X_1},0_{Y_1})$. Thus $F'_3(id_{\mathcal{H}_1})=id_{F_3(\mathcal{H}_1)}$. Therefore $F'_3$ is a functor.
\begin{pro}\label{16}
The categories ${\bf LFtrans}^\downarrow$ and ${\bf LSpaceFP}$ are isomorphic.
\end{pro}
\textbf{Proof:} Let $F_3:{{\bf LSpaceFP}}\rightarrow {{\bf LFtrans}^{\downarrow}}$ and $F'_3:{{\bf LFtrans}^{\downarrow}}\rightarrow {{\bf LSpaceFP}}$ be functors. Then $F'_3\circ F_3:{{\bf LSpaceFP}}\rightarrow {\bf LSpaceFP}$ is a functor such that for all $\pi\in|{\bf LSpaceFP}|$ and ${\bf LSpaceFP}$-morphism $(\phi,\mu):\pi_1\rightarrow\pi_2$, \begin{eqnarray*}
(F'_3\circ F_3)(\pi_1)&=&F'_3( F_3(\pi_1))\\
&=&F'_3(\mathcal{H}_{\pi_1})\\
&=&\pi_{\mathcal{H}_{\pi_1}}\\
&=&\pi_1.
\end{eqnarray*}
Thus $(F'_3\circ F_3)(\pi_1)=\pi_1$ and 
\begin{eqnarray*}
(F'_3\circ F_3)(\phi,\mu)&=&F'_3( F_3(\phi,\mu))\\
&=&F'_3(\phi,\mu)\\
&=&(\phi,\mu).
\end{eqnarray*}
Thus $(F'_3\circ F_3)(\phi,\mu)=(\phi,\mu)$ and therefore $F'_3\circ F_3=id_{{\bf LSpaceFP}}$. Similarly, we can show that $F_3\circ F'_3:{{\bf LFtrans}^{\downarrow}}\rightarrow {{\bf LFtrans}^{\downarrow}}$ is a functor such that for all $\mathcal{H}_1\in|{{\bf LFtrans}^{\downarrow}}|$ and ${{\bf LFtrans}^{\downarrow}}$-morphism $(\psi,\nu):\mathcal{H}_1\rightarrow\mathcal{H}_2$, $(F_3\circ F'_3)(\mathcal{H}_1)=\mathcal{H}_1,(F_3\circ F'_3)(\psi,\nu)=(\psi,\nu)$.
Thus $F_3\circ F'_3=id_{{\bf LFtrans}^{\downarrow}}.$
Hence the category ${\bf LSpaceFP}$ is isomorphic to the category ${\bf LFtrans}^\downarrow$.\\\\
{The connections between the categories ${\bf LFtrans}^\downarrow$ and ${\bf Ftrans}^\downarrow$ with the category ${\bf Qua}$ are described in the following.}
\begin{pro}\label{17} There is a subcategory of the category ${\bf Qua}$, which is isomorphic to  the category ${\bf LFtrans}^{\downarrow}$ .
\end{pro}
\textbf{Proof:} {Propositions \ref{5} and \ref{16} lead to this proof.}
\begin{pro}\label{18} The subcategory of a category, which is isomorphic to a subcategory of the category ${\bf Qua}$, is isomorphic to the category ${\bf Ftrans}^{\downarrow}$.
\end{pro}
\textbf{Proof:} {Propositions \ref{5}, \ref{0} and \ref{17} lead to this proof.}
\section{Categorical view of $L$-valued fuzzy partition with $L$-valued fuzzy pretopology and $L$-valued fuzzy interior operator}
In this section, we introduce the concepts of the categories of $L$-valued fuzzy pretopological spaces with morphisms as pairs of $L$-valued fuzzy relations, as well as the category of \v{C}ech $L$-valued fuzzy interior spaces. Also, we demonstrate the relationships of the above-introduced categories with ${\bf Qua}$. Now, we recall the following from \cite{qi,zh}.
\begin{def1}\label{PTOP}
An {\bf $L$-valued fuzzy pretopology} on a nonempty set $X$ is a set of functions $P=\lbrace p_x:L^X\rightarrow L|x\in X\rbrace$ such that 
\item[(i)] $p_x(\textbf{a})= a,\,\forall\,a\in L,\textbf{a}\in L^X$, 
\item[(ii)] $p_x(f)\leq f(x),\forall\,f\in L^X$, and
\item[(iii)] $p_x(f\wedge g)=p_x(f)\wedge p_x(g),\forall\,f,g\in L^X$.
\\\\
The pair $\Sigma=(X,P)$ is called an {\bf $L$-valued fuzzy pretopological space}.
\end{def1}
\begin{def1}
Let $\Sigma_1=(X_1,{P}_1)$ and $\Sigma_2=(X_2,{P}_2)$ be $L$-valued fuzzy pretopological spaces, where ${{P}_1}=\{p_{x_1}:L^{X_1}\rightarrow L|x_1\in X_1\}$ and ${{P}_2}=\{p_{x_2}:L^{X_2}\rightarrow L:x_2\in X_2\}$ are $L$-valued fuzzy pretopologies on $X_1$ and $X_2$, respectively. Then an {\bf $L$-valued fuzzy continuous function} from $\Sigma_1$ to $\Sigma_2$ is a pair $(\varrho,\varrho')$ such that
\begin{itemize}
\item[(i)] $\varrho:X_1\times X_2\rightarrow L$ and $\varrho':X_2\times X_1\rightarrow L$ are $L$-valued fuzzy relations, and
\item[(ii)] $\neg\varrho'(x_2,x_1)\ast p_{x_2}(f)\leq  p_{x_1}(\overleftarrow{\varrho}(f)),\,\forall\,x_1\in X_1,x_2\in X_2,f\in L^{X_2}$.
\end{itemize}
\end{def1}
Let $\Sigma_1=(X_1,{P}_1),\Sigma_2=(X_2,{P}_2),\Sigma_3=(X_3,{P}{_3})$ be $L$-valued fuzzy pretopological spaces and $(\varrho_1,\varrho'_1):\Sigma_1\rightarrow \Sigma_2,(\varrho_2,\varrho'_2):\Sigma_2\rightarrow \Sigma_3$ be their $L$-valued fuzzy continuous functions. Then the {\bf composition}  $(\varrho_2,\varrho'_2)\bullet  (\varrho_1,\varrho'_1):\Sigma\rightarrow \Sigma_3$ of $(\varrho_1,\varrho'_1)$ and $(\varrho_2,\varrho'_2) $ is a pair of $L$-valued fuzzy relations such that $ (\varrho_2,\varrho'_2)\bullet  (\varrho_1,\varrho'_1)=(\varrho_2\odot \varrho_1,\varrho'_1\oplus\varrho'_2)$.
\begin{pro}\label{19} $L$-valued fuzzy pretopological spaces with their $L$-valued fuzzy continuous functions form a category.
\end{pro}
\textbf{Proof:} It is enough to demonstrate that the composition of $L$-valued fuzzy continuous functions is also an $L$-valued fuzzy continuous function. To do this, let $\Sigma_1=(X_1,P_1),\Sigma_2=(X_2,P{_2}),\Sigma_3=(X_3,P{_3}) $ be $L$-valued fuzzy pretopological spaces and $(\varrho_1,\varrho'_1):\Sigma_1\rightarrow \Sigma_2 $, $(\varrho_2,\varrho'_2):\Sigma_2\rightarrow \Sigma_3$ be $L$-valued fuzzy continuous functions. Then for all $x_1\in X_1,f\in L^{X_3}$
\begin{eqnarray*}
p_{x_1}(\overleftarrow{\varrho_2\odot \varrho_1}(f))&=&p_{x_1}((\overleftarrow{\varrho_1}\circ \overleftarrow{\varrho_2})(f))\\
&=&p_{x_1}(\overleftarrow{\varrho_1}( \overleftarrow{\varrho_2}(f)))\\
&\geq&\neg\varrho'_1(x_2,x_1)\ast {p_{x_2}}( \overleftarrow{\varrho_2}(f))\\
&\geq&\neg\varrho'_1(x_2,x_1)\ast \neg\varrho'_2(x_3,x_2)\ast p_{x_3}( f)\\
&=&\neg(\varrho'_1(x_2,x_1)\#\varrho'_2(x_3,x_2))\ast p_{x_3}( f)\\
&\geq&\bigvee\limits_{x_2\in X_2}\neg(\varrho'_1(x_2,x_1)\#\varrho'_2(x_3,x_2))\ast p_{x_3}( f)\\
&=&\neg(\bigwedge\limits_{x_2\in X_2}\varrho'_1(x_2,x_1)\#\varrho'_2(x_3,x_2))\ast p_{x_3}( f)\\
&=&\neg(\varrho'_1\oplus\varrho'_2)(x_3,x_1)\ast p_{x_3}( f).
\end{eqnarray*}
Thus $(\varrho_2\odot \varrho_1,\varrho'_1\oplus\varrho'_2):\Sigma_1\rightarrow\Sigma_3$ is an $L$-valued fuzzy continuous function.\\\\
We shall denote by ${\bf LFPrTop}$, the category of $L$-valued fuzzy pretopological spaces with $L$-valued fuzzy continuous function as morphisms.\\\\
In the following, we show that there exists a subcategory of the category ${\bf LFPrTop}$, which is isomorphic to the category of $L$-valued fuzzy pretopological spaces (say, ${\bf FPrTop}$) introduced in \cite{mockor}, whose morphisms are pairs of functions between the underlying sets of corresponding objects. For which, we assume that ${\bf SLFPrTop}$ is a category, whose objects are same as the objects of ${\bf LFPrTop}$ and morphisms are ${\bf LFPrTop}$-morphisms $(\varrho,\varrho'):\Sigma_1\rightarrow \Sigma_2$ such that
\begin{itemize}
\item[$\bullet$] $\varrho:X_1\times X_2\rightarrow L$ is an $L$-valued fuzzy relation and for all $x_1\in X_1$, there exists a unique $x_2\in X_2, \varrho(x_1,x_2)=1$ and 0, otherwise,
\item[$\bullet$] $\varrho':X_2\times X_1\rightarrow L$ is an $L$-valued fuzzy relation and for all $x_2\in X_2$ and $x_1\in X_1,\,\varrho'(x_2,x_1)\in\{0,1\}$.
\end{itemize}
{It is obvious that ${\bf SLFPrTop}$ is a subcategory of} ${\bf LFPrTop}$.
\begin{pro}\label{20} Let $F_4:{{\bf FPrTop}}\rightarrow {{\bf SLFPrTop}}$ be a function such that for every $\Sigma_1\in |{\bf FPrTop}|,F_4(\Sigma_1)=\Sigma_1$ and for every ${\bf FPrTop}$-morphism $\zeta:\Sigma_1\rightarrow \Sigma_2,F_4(\zeta):\Sigma_1\rightarrow \Sigma_2$ be a function such that $F_4(\zeta)=(\varrho_{\zeta},\varrho'_{\zeta})$, where $\varrho_{\zeta}:X_1\times X_2\rightarrow L,\varrho'_{\zeta}:X_2\times X_1\rightarrow L$ are $L$-valued fuzzy relations such that for each $x_1\in X_1,x_2\in X_2$,
\begin{eqnarray*} 
\varrho_{\zeta}(x_1,x_2)=\begin{cases}
1 &\text{ if } \zeta(x_1)=x_2\\
0&\text{ otherwise,}
\end{cases} 
 \varrho'_{\zeta}(x_2,x_1)=\begin{cases}
0 &\text{ if } \zeta^{-1}(x_2)=x_1\\
1&\text{ otherwise,}
\end{cases}
\end{eqnarray*}
respectively. Then $F_4$ is a functor.
\end{pro}
\textbf{Proof:} {It is easy to understand the proof.}
\begin{pro}\label{21} Let $F'_4:{{\bf SLFPrTop}}\rightarrow {{\bf FPrTop}}$ be a function such that for every $\Sigma_1\in |{\bf SLFPrTop}|,F'_4(\Sigma_1)=\Sigma_1$ and for every ${\bf SLFPrTop}$-morphism $(\varrho,\varrho'):\Sigma_1\rightarrow \Sigma_2,F'_4(\varrho,\varrho'):\Sigma_1\rightarrow \Sigma_2$ be a function such that $F_4'(\varrho,\varrho')={\zeta}_{\varrho},$ where $\zeta_{\varrho}:X_1\rightarrow X_2$ is a function such that for each $x_1\in X_1,\zeta_{\varrho}(x_1)=x_2$ iff $\varrho(x_1,x_2)=1$. Then $F_4'$ is a functor. 
\end{pro}
\textbf{Proof:} {It is easy to understand the proof.}
\begin{pro}\label{22} The categories ${\bf FPrTop}$ and  ${\bf SLFPrTop}$ are isomorphic.
\end{pro}
\textbf{Proof:} Let $F_4:{{\bf FPrTop}}\rightarrow {{\bf SLFPrTop}}$ and $F'_4:{\bf SLFPrTop}\rightarrow {{\bf FPrTop}}$ be functors. Then $F'_4\circ F_4:{{\bf FPrTop}}\rightarrow {{\bf FPrTop}}$ is a functor such that for all $\Sigma_1\in|{{\bf FPrTop}}|$ and ${{\bf FPrTop}}$-morphism $\zeta:\Sigma_1\rightarrow\Sigma_2$, 
\begin{eqnarray*}
(F'_4\circ F_4)(\Sigma_1)&=&
F'_4( F_4(\Sigma_1))\\
&=&F'_4(\Sigma_1)\\
&=&\Sigma_1.
\end{eqnarray*}
Thus $(F'_4\circ F_4)(\Sigma_1)=\Sigma_1$ and 
\begin{eqnarray*}
(F'_4\circ F_4)(\zeta)&=&F'_4( F_4(\zeta))\\
&=&F'_4(\varrho_\zeta,\varrho'_\zeta)\\
&=&\zeta_{\varrho_\zeta}\\
&=&\zeta.
\end{eqnarray*}
Thus $(F'_4\circ F_4)(\zeta)=\zeta$ and therefore $F'_4\circ F_4=id_{{\bf FPrTop}}$. Similarly, we can show that $F_4\circ F'_4:{{\bf SLFPrTop}}\rightarrow {\bf SLFPrTop}$ is a functor such that for all $\Sigma_1\in|{{\bf SLFPrTop}}|$ and ${{\bf SLFPrTop}}$-morphism $(\varrho,\varrho'):\Sigma_1\rightarrow\Sigma_2,(F_4\circ F'_4)(\Sigma_1)=\Sigma_1,(F_4\circ F'_4)(\varrho,\varrho')=(\varrho,\varrho')$.
Thus $F_4\circ F'_4=id_{{\bf SLFPrTop}}.$
Therefore the categories ${\bf FPrTop}$ and ${\bf SLFPrTop}$ are isomorphic.\\\\
In the following, we unify categories ${\bf LFPrTop}$ and ${\bf FPrTop}$ as the ${\bf Qua}$ category.
\begin{pro}\label{25} The category ${\bf LFPrTop}$ is a subcategory of the category ${\bf Qua}$.
\end{pro}
\textbf{Proof:} Let $\Sigma_1=(X_1,{P}_1)\in |{\bf LFPrTop}|$ and $(\varrho,\varrho'):\Sigma_1\rightarrow \Sigma_2$ be an ${\bf LFPrTop}$-morphism. Then $\Sigma$ can be seen as an object $ \mathcal{M}_{\Sigma_1}=(M^{*},M_{*},\mathbf{M}_{P_1})\in|{\bf Qua}|$, where $M^{*}=X_1=M_{*}$ and $\mathbf{M}_{P_1}:X_1\times X_1\rightarrow L$ is an $L$-valued fuzzy relation such that $\mathbf{M}_{P_1}(x_1,x_{1}')=\neg p_{x_1}(\neg1_{\{x_1'\}}),\,\forall\,x_1,x_1'\in X$. Also, we can construct a pair of functions $(g_{*},g^{*}):\mathcal{M}_{\Sigma_1}\rightarrow \mathcal{M}_{\Sigma_2}$, where $ g_{*}:X_1\times X_2\rightarrow L,g^{*}:X_2\times X_1\rightarrow L$ are $L$-valued fuzzy relations such that $g_{*}(x_1,x_2)=\varrho(x_1,x_2),g^{*}(x_2,x_1)=\varrho'(x_2,x_1),\,\forall\,x_1\in X_1,x_2\in X_2$. Now, let $x_1\in X_1,x_2\in X_2$. Then
\begin{eqnarray*}
p_{x_1}(\overleftarrow{\varrho}(\neg1_{\{x_2\}}))&=&p_{x_1}(\bigwedge\limits_{x'_{1}\in X}\overleftarrow{\varrho}(\neg1_{\{x_2\}})(x_1')\#\neg 1_{\{x_1'\}})\\
&=& \bigwedge\limits_{x_1'\in X}\overleftarrow{\varrho}(\neg1_{\{x_2\}})(x_1')\#p_{x_1}(\neg 1_{\{x_1'\}})\\
&=&\bigwedge\limits_{x_1'\in X}\bigwedge\limits_{x'_2\in X_2}(\neg\varrho(x_1',x'_2)\#\neg1_{\{x_2\}}(x'_2))\#p_{x_1}(\neg 1_{\{x_1'\}})\\
&=&\bigwedge\limits_{x_1'\in X}\neg\varrho(x_1',x_2)\#p_{x_1}(\neg 1_{\{x_1'\}})\\
&\leq&\neg\varrho(x_1',x_2)\#p_{x_1}(\neg 1_{\{x_1'\}})\\
&=&\neg g_{*}(x_1',x_2)\#\neg\mathbf{M}_{P_1}(x_1,x_1')\\
&=&\neg (g_{*}(x_1',x_2)\ast \mathbf{M}_{P_1}(x_1,x_1')).
\end{eqnarray*}
On the other hand,
\begin{eqnarray*}
g^{*}(x_2',x_1)\# \mathbf{M}_{P_2}(x_2',x_2)&=&\varrho'(x_2',x_1)\#\neg p_{x_2'}(\neg1_{\{x_2\}})\\
&=&\neg(\neg\varrho'(x_2',x_1)\ast p_{x_2'}(\neg1_{\{x_2\}}))\\
&\geq&\neg(p_{x_1}(\overleftarrow{\varrho}(\neg1_{\{x_1'\}})))\\
&\geq&\neg\neg (g_{*}(x_1',x_2)\ast \mathbf{M}_P(x_1,x_1')).\\
&=&g_{*}(x_1',x_2)\ast \mathbf{M}_P(x_1,x_1').
\end{eqnarray*}
Thus $ g_{*}(x_1',x_2)\ast \mathbf{M}_P(x_1,x_1')\leq g^{*}(x_2',x_1)\# \mathbf{M}_{P'}(x_2',x_2)$.
Therefore $(g_{*},g^{*}):\mathcal{M}_{\Sigma_1}\rightarrow \mathcal{M}_{\Sigma_2}$ is a ${\bf Qua}$-morphism. Hence the category ${\bf LFPrTop}$ is a subcategory of the category ${\bf Qua}$.
\begin{pro}\label{26} There is a subcategory of the category ${\bf Qua}$, which is isomorphic to the category ${\bf FPrTop}$.
\end{pro}
\textbf{Proof:} {Propositions \ref{22} and \ref{25} lead to this proof.}\\\\
Next, we recall the following from \cite{mockor,spt1}.
\begin{def1}\label{CINT}
A function ${i}:L^X\rightarrow L^X$ is called {\bf \v{C}ech $L$-valued fuzzy interior operator} if
\item[(i)] ${i}(\textbf{a})=\textbf{a},\,\forall \,\textbf{a}\in L^X$,
\item[(ii)] ${i}(f)\leq f,\,\forall\,f\in L^X$, and
\item[(iii)] ${i}(f\wedge g)=i(f)\wedge {i}(g),\,\forall\,f,g\in L^X.$\\\\
The pair $I=(X,i)$ is called {\bf \v{C}ech $L$-valued fuzzy interior space}.
\end{def1}
{In the following, the concept of the category of \v{C}ech $L$-valued fuzzy interior spaces is introduced and their morphisms are defined as pairs of $L$-valued fuzzy relations between the underlying sets of corresponding objects.} We start with the following.
\begin{def1}
Let $I_1=(X_1,i_1)$ and $I_2=(X_2,i{_2})$ be \v{C}ech $L$-valued fuzzy interior spaces, where ${i}_1:L^{X_1}\rightarrow L^{X_1}$ and ${i{_2}}:L^{X_2}\rightarrow L^{X_2}$ are \v{C}ech $L$-valued fuzzy interior operators on $X_1$ and $X_2$, respectively. Then an {\bf $L$-valued fuzzy continuous function} from ${I}_1$ to ${I}_2$ is a pair $(\kappa,\kappa')$ such that
\begin{itemize}
\item[(i)] $\kappa:X_1\times X_2\rightarrow L$ and $\kappa':X_2\times X_1\rightarrow L$ are $L$-valued fuzzy relations, and
\item[(ii)] $\neg\kappa'(x_2,x_1)\ast i{_2}(f)(x_2)\leq  i_1(\overleftarrow{\kappa}(f))(x_1),\,\forall\,x_1\in X_1,x_2\in X_2,f\in L^{X_2}$.
\end{itemize} 
\end{def1}
Let ${I}_1=(X_1,i_1),{I}_2=(X_2,i_2),{I}_3=(X_3,{i}{_3})$ be \v{C}ech $L$-valued fuzzy interior spaces and $(\kappa_1,\kappa'_1):{I}_1\rightarrow {I}_2,(\kappa_2,\kappa'_2):{I}_2\rightarrow {I}_3$ be their $L$-valued fuzzy continuous functions. Then the {\bf composition} $(\kappa_2,\kappa'_2)\bullet  (\kappa_1,\kappa'_1):{I}_1\rightarrow {I}_3$ of $ (\kappa_1,\kappa'_1)$ and $(\kappa_2,\kappa'_2) $ is a pair of $L$-valued fuzzy relations such that $ (\kappa_2,\kappa'_2)\bullet  (\kappa_1,\kappa'_1)=(\kappa_2\odot \kappa_1,\kappa'_1\oplus\kappa'_2)$.
\begin{pro} \v{C}ech $L$-valued fuzzy interior spaces with their $L$-valued fuzzy continuous functions form a category.
\end{pro}
\textbf{Proof:} It is enough to demonstrate that composition of $L$-valued fuzzy continuous functions is also an $L$-valued fuzzy continuous function. To do this, let ${I_1}=(X_1,i_1),{I}_2=(X_2,i{_2}),{I}_3=(X_3,i{_3}) $ be \v{C}ech $L$-valued fuzzy interior spaces and $(\kappa_1,\kappa'_1):{I}_1\rightarrow {I}_2 $, $(\kappa_2,\kappa'_2):{I}_2\rightarrow {I}_3$ be $L$-valued fuzzy continuous functions. Then for all $x_1\in X_1,f\in L^{X3}$
\begin{eqnarray*}
i_1(\overleftarrow{\kappa_2\odot \kappa_1}(f))(x_1)&=&i_1((\overleftarrow{\kappa_1}\circ \overleftarrow{\kappa_2})(f))(x_1)\\
&=&i_1(\overleftarrow{\kappa_1}( \overleftarrow{\kappa_2}(f)))(x_1)\\
&\geq&\neg\kappa'_1(x_2,x_1)\ast {i{_2}}( \overleftarrow{\kappa_2}(f))(x_2)\\
&\geq&\neg\kappa'_1(x_2,x_1)\ast \neg\kappa'_2(x_3,x_2)\ast i{_3}( f)(x_3)\\
&=&\neg(\kappa'_1(x_2,x_1)\#\kappa'_2(x_3,x_2))\ast i{_3}( f)(x_3)\\
&\geq&\bigvee\limits_{x_2\in X_2}\neg(\kappa'_1(x_2,x_1)\#\kappa'_2(x_3,x_2))\ast i{_3}( f)(x_3)\\
&=&\neg(\bigwedge\limits_{x_2\in X_2}\kappa'_1(x_2,x_1)\#\kappa'_2(x_3,x_2))\ast i{_3}( f)(x_3)\\
&=&\neg(\kappa'_1\oplus\kappa'_2)(x_3,x_1)\ast i{_3}( f)(x_3).
\end{eqnarray*}
Thus $(\kappa_2\odot \kappa_1,\kappa'_1\oplus\kappa'_2):{I}_1\rightarrow{I}_3$ is an $L$-valued fuzzy continuous function.\\\\
We shall denote by ${\bf LFCInt}$, the category of \v{C}ech $L$-valued fuzzy interior spaces with $L$-valued fuzzy continuous functions as morphisms.\\\\
In the following, we show that there exists a subcategory of the category ${\bf LFCInt}$, which is isomorphic to the category of \v{C}ech $L$-valued fuzzy interior spaces (say, ${\bf CInt}$) introduced in \cite{mockor}, whose morphisms are given as functions between the underlying sets of corresponding objects. For which, we assume that ${\bf SLFCInt}$ is a category,  whose objects are same as the objects of ${\bf LFCInt}$ and morphisms are ${\bf LFCInt}$-morphisms $(\kappa,\kappa'):{I}_1\rightarrow {I}_2$ such that
\begin{itemize}
\item[$\bullet$] $\kappa:X_1\times X_2\rightarrow L$ is an $L$-valued fuzzy relation such that for each $x_1\in X_1$ there exists a unique $x_2\in X_2, \kappa(x_1,x_2)=1$ and 0, otherwise,
\item[$\bullet$] $\kappa':X_2\times X_1\rightarrow L$ is an $L$-valued fuzzy relation such that for each $x_2\in X_2$ and $x_1\in X_1,\,\kappa'(x_2,x_1)\in\{0,1\}$.
\end{itemize}
{It is obvious that ${\bf SLFCInt}$ is a subcategory of} ${\bf LFCInt}$.
\begin{pro} Let $F_5:{{\bf CInt}}\rightarrow {{\bf SLFCInt}}$ be a function such that for every ${I}_1\in |{\bf CInt}|,F_5({I}_1)={I}_1$ and for every ${\bf CInt}$-morphism $\omega:{I}_1\rightarrow {I}_2,F_5(\omega):{I}_1\rightarrow {I}_2$ be a function such that $F_5(\omega)=(\kappa_{\omega},\kappa'_{\omega})$, where $\kappa_{\omega}:X_1\times X_2\rightarrow L,\kappa'_{\omega}:X_2\times X_1\rightarrow L$ are $L$-valued fuzzy relations such that for each $x_1\in X_1,x_2\in X_2$,
\begin{eqnarray*} 
\kappa_{\omega}(x_1,x_2)=\begin{cases}
1 &\text{ if } \omega(x_1)=x_2\\
0&\text{ otherwise,} 
\end{cases}
 \kappa'_{\omega}(x_2,x_1)=\begin{cases}
0 &\text{ if } \omega^{-1}(x_2)=x_1\\
1&\text{ otherwise,}
\end{cases}
\end{eqnarray*}
respectively. Then $F_5$ is a functor.
\end{pro}
\textbf{Proof:} {It is easy to understand the proof.}
\begin{pro} Let $F'_5:{{\bf SLFCInt}}\rightarrow {{\bf CInt}}$ be a function such that for every ${I}_1\in |{\bf SLFCInt}|,F'_5({I_1})={I_1}$ and for every ${\bf SLFCInt}$-morphism $(\kappa,\kappa'):{I_1}\rightarrow {I}_2,F'_5(\kappa,\kappa'):{I}_1\rightarrow {I}_2$ be a function such that $F_5'(\kappa,\kappa')={\omega}_{\kappa}$, where $\omega_{\kappa}:X_1\rightarrow X_2$ is a function such that for each $x_1\in X_1,\omega_{\kappa}(x_1)=x_2$ iff $\kappa(x_1,x_2)=1$. Then $F_5'$ is a functor. 
\end{pro}
\textbf{Proof:} {It is easy to understand the proof.}
\begin{pro}\label{isoclc} The categories ${\bf CInt}$ and ${\bf SLFCInt}$ are isomorphic.
\end{pro}
\textbf{Proof:} Let $F_5:{{\bf CInt}}\rightarrow {{\bf SLFCInt}}$ and $F'_5:{\bf SLFCInt}\rightarrow {{\bf CInt}}$ be functors. Then $F'_5\circ F_5:{{\bf CInt}}\rightarrow {{\bf CInt}}$ is a functor such that for all $I_1\in|{{\bf CInt}}|$ and ${{\bf CInt}}$-morphism $\omega:I_1\rightarrow I_2$, 
\begin{eqnarray*}
(F'_5\circ F_5)(I_1)&=&F'_5( F_5(I_1))\\
&=&F'_5(I_1)\\
&=&I_1.
\end{eqnarray*}
Thus $(F'_5\circ F_5)(I_1)=I_1$ and 
\begin{eqnarray*}
(F'_5\circ F_5)(\omega)&=&F'_5( F_5(\zeta))\\
&=&F'_5(\kappa_\omega,\kappa'_\omega)\\
&=&\omega_{\kappa_\omega}\\
&=&\omega.
\end{eqnarray*}
Thus $(F'_5\circ F_5)(\omega)=\omega$ and therefore $F'_5\circ F_5=id_{{\bf CInt}}$. Similarly, we can show that $F_5\circ F'_5:{{\bf SLFCInt}}\rightarrow {\bf SLFCInt}$ is a functor such that for all $I_1\in|{{\bf SLFCInt}}|$ and ${{\bf SLFCInt}}$-morphism $(\kappa,\kappa'):I_1\rightarrow I_2,(F_5\circ F'_5)(I_1)=I_1,(F_5\circ F'_5)(\kappa,\kappa')=(\kappa,\kappa')$.
Thus $F_5\circ F'_5=id_{{\bf SLFCInt}}.$
therefore categories ${\bf CInt}$ and ${\bf SLFCInt}$ are isomorphic.\\\\
Next, we show establish an isomorphism  between the category ${\bf LFPrTop}$ and the category ${\bf LFCInt}$. To do this, let $\Sigma=(X,{P})$ be an $L$-valued fuzzy pretopological space and $I=(X,i)$ be a \v{C}ech $L$-valued fuzzy interior space. Then we can introduce a \v{C}ech $L$-valued fuzzy interior space and $L$-valued fuzzy pretopological space with respect to $\Sigma$ and $I$, respectively, defined by $I_\Sigma$ and $\Sigma_I$, as follows:
 \begin{itemize}
     \item[(i)] $I_\Sigma=(X,i_{P})$, where $i_{P}:L^X\rightarrow L^X$ such that for each $x\in X,f\in L^X$, $i_{{P}}(f)(x)=p_{x_1}(f),$ and
  \item[(ii)] $\Sigma_{I}=(X,{P}_{i})$, where ${P}_{i}=\{p^{i}_{x}:L^X\rightarrow L|x\in X\}$ such that for each $x\in X,f\in L^X,\,p^{i}_x(f)=i(f)(x)$.
 \end{itemize}
 \begin{pro}
Each ${\bf LFPrTop}$-morphism is also an ${\bf LFCInt}$-morphism.
\end{pro}
\textbf{Proof:} Let $\Sigma_1=(X_1,P_1),\Sigma_2=(X_2,P{_2})\in|{\bf LFPrTop}|$ and $(\varrho,\varrho'):\Sigma_1\rightarrow \Sigma_2$ be an ${\bf LFPrTop}$-morphism. Then for all $x_1\in X_1,f\in L^{X_2}$
\begin{eqnarray*}
i_{P_1}(\overleftarrow{\varrho}(f))(x_1)&=& p_{x_1}(\overleftarrow{\varrho}(f))\\
&\geq& \neg\varrho'(x_2,x_1)\ast p_{x_2}(f)\\
&=& \neg\varrho'(x_2,x_1)\ast i_{P{{_2}}}(f)(x_2).
\end{eqnarray*}
Thus $i_{P_1}(\overleftarrow{\varrho}(f))(x_1)\geq \neg\varrho'(x_2,x_1)\ast i_{P{_2}}(f)(x_2)$. Therefore $(\varrho,\varrho'):I_{\Sigma_1}\rightarrow I_{\Sigma_2}$ is an ${\bf LFCInt}$-morphism.
\begin{pro}
Each ${\bf LFCInt}$-morphism is also an ${\bf LFPrTop}$-morphism.
\end{pro}
\textbf{Proof:} Let $I_1=(X_1,i_1),I_2=(X_2,i{_2})\in|{\bf LFCInt}|$ and $(\kappa,\kappa'):I_1\rightarrow I_2$ be an ${\bf LFCInt}$-morphism. Then for all $x_1\in X_1,f\in L^{X_2}$
\begin{eqnarray*}
p^{i_1}_{x_1}(\overleftarrow{\kappa}(f))&=&i_1(\overleftarrow{\kappa}(f))(x_1)\\
&\geq&\neg\kappa'(x_2,x_1)\ast i{{_2}}(f)(x_2)\\
&=& \neg\kappa'(x_2,x_1)\ast p^{i{_2}}_{x_2}(f).
\end{eqnarray*}
Thus $p^{i_1}_{x_1}(\overleftarrow{\kappa}(f))\geq \neg\kappa'(x_2,x_1)\ast p^{i{_2}}_{x_2}(f)$. Therefore $(\kappa,\kappa'):\Sigma_{I_1}\rightarrow {\Sigma_{I_2}}$ is an ${\bf LFCInt}$-morphism.
\begin{pro}\label{pi}
Let $F_6:{\bf LFPrTop}\rightarrow {\bf LFCInt}$ be a function such that for every $\Sigma_1\in |{\bf LFPrTop}|, F_6(\Sigma_1)=I_{\Sigma_1}$ and for every ${\bf LFPrTop}$-morphism $(\varrho,\varrho'):\Sigma_1\rightarrow\Sigma_2,F_6(\varrho,\varrho'):I_{\Sigma_1}\rightarrow I_{\Sigma_2}$ be a function such that $F_6(\varrho,\varrho')=(\varrho,\varrho')$. Then $F_6$ is a functor.
\end{pro}
\textbf{Proof:} (i) Let $(\varrho_1,\varrho'_1):\Sigma_1\rightarrow\Sigma_2$ and $(\varrho_2,\varrho'_2):\Sigma_2\rightarrow\Sigma_3$ be  ${\bf LFPrTop}$-morphisms. Then 
\begin{eqnarray*}
F_6((\varrho_2,\varrho'_2)\bullet(\varrho_1,\varrho'_1))&=&(\varrho_2,\varrho'_2)\bullet(\varrho_1,\varrho'_1)\\
&=&F_6(\varrho_2,\varrho'_2)\bullet F_6(\varrho_1,\varrho'_1).
\end{eqnarray*}
Therefore $F_6((\varrho_2,\varrho'_2)\bullet(\varrho_1,\varrho'_1))=F_6(\varrho_2,\varrho'_2)\bullet F_6(\varrho_1,\varrho'_1)$.\\\\
(ii) Let $\Sigma_1=(X_1,P_1)\in|{\bf LFPrTop}|$. Then $F_6(id_{\Sigma_1})=F_6(1_{X_1},0_{X_1})=(1_{X_1},0_{X_1})$. Thus $F_6(id_{\Sigma_1})=id_{F_6({\Sigma_1})}$. Hence $F_6$ is a functor.
\begin{pro}
Let $F'_6:{\bf LFCInt}\rightarrow {\bf LFPrTop}$ be a function such that for every $I\in |{\bf LFCInt}|, F'_6(I_1)={\Sigma}_{I_1}$ and for every ${\bf LFCInt}$-morphism $(\kappa,\kappa'):I_1\rightarrow I_2,F'_6(\kappa,\kappa'):{\Sigma_{I_1}}\rightarrow {\Sigma_{I_2}}$ be a function such that $F'_6(\kappa,\kappa')=(\kappa,\kappa')$. Then $F'_6$ is a functor.
\end{pro}
\textbf{Proof:} Similar to that of Proposition \ref{pi}.
\begin{pro}\label{isopi}
The category ${\bf LFPrTop}$ and ${\bf LFCInt}$ are ismorphic.
\end{pro}
\textbf{Proof:} Let $F_6:{\bf LFPrTop}\rightarrow {\bf LFCInt}$ and $F'_6:{\bf LFCInt}\rightarrow {\bf LFPrTop}$ be functors. Then $F'_6\bullet F_6:{\bf LFPrTop}\rightarrow {\bf LFPrTop}$ is a functor such that for all $\Sigma_1\in |{\bf LFPrTop}|$ and ${\bf LFPrTop}$-morphism $(\varrho,\varrho'):\Sigma_1\rightarrow \Sigma_2$,
\begin{eqnarray*}
(F'_6\bullet F_6)(\Sigma_1)&=&F'_6(F_6(\Sigma_1))\\
&=&F'_6(I_{\Sigma_1})\\
&=&\Sigma_{I_{\Sigma_1}}\\
&=&\Sigma_1.
\end{eqnarray*}
Thus $(F'_6\bullet F_6)(\Sigma_1)=\Sigma_1\,$ and 
\begin{eqnarray*}
(F'_6\bullet F_6)(\varrho,\varrho')&=&F'_6(F_6(\varrho,\varrho'))\\
&=&F'_6(\varrho,\varrho')\\
&=&(\varrho,\varrho').
\end{eqnarray*}
Therefore $F'_6\bullet F_6=id_{\bf LFPrTop}$. Similarly, we can show that $F_6\bullet F'_6:{\bf LFCInt}\rightarrow {\bf LFCInt}$ is a functor such that $F_6\bullet F'_6={id_{\bf LFCInt}}$. Therefore the categories ${\bf LFPrTop}$ and ${\bf LFCInt}$ are isomorphic.\\\\
{The connections between the categories ${\bf LFCInt}$ and ${\bf CInt}$ with the category ${\bf Qua}$ are described in the following.}
\begin{pro}\label{isop2} There is a subcategory of the category ${\bf Qua}$, which is isomorphic to the category ${\bf LFCInt}$.
\end{pro}
\textbf{Proof:} {Propositions \ref{25} and \ref{isopi} lead to this proof.}
\begin{pro} The subcategory of a category, which is isomorphic to a subcategory of the category ${\bf Qua}$, is isomorphic to the category ${\bf CInt}$.
\end{pro}
\textbf{Proof:} { Propositions \ref{25}, \ref{isoclc} and \ref{isop2} lead to this proof.}\\\\
For the following, let $\pi = (X,J,\mathcal{P})$ be a space with $L$-valued fuzzy partition,  $\mathcal{H}=(X,J,v,H)$ be an $L$-valued fuzzy lower transformation system, $\Sigma=(X,P) $ be an $L$-valued fuzzy pretopological space and $I=(X,i)$ be a \v{C}ech $L$-valued fuzzy interior space. Then we can introduce an $L$-valued fuzzy pretopological space and a \v{C}ech $L$-valued fuzzy interior space corresponding to $\pi$ (respectively, $\mathcal{H}$), denoted by $\Sigma_\pi$ and $I_\pi$ (respectively, $\Sigma_\mathcal{H}$ and $I_\mathcal{H}$), as follows: 
\begin{itemize}
\item[(i)] ${\Sigma}_{\pi}=(X, P_{\mathcal{P}})\in|{\bf LFPrTop}|$, where ${P}_{\mathcal{P}}=\{p^{{\mathcal{P}}}_x:L^X\rightarrow L|x\in X\}$ is an $L$-valued fuzzy pretopology such that for all $x\in X,f\in L^X,$
\[p^{{\mathcal{P}}}_{x}(f)=\bigwedge\limits_{x'\in X}\neg A_{\xi(x)}(x')\#f(x')=F_{\mathcal{P}}^{\downarrow}[f](\xi(x)),\]
where $\xi(x)\in J$ is the unique index such that $x\in core(A_{\xi(x)}),$ and
\item[(ii)] ${I}_{\pi}=(X, P_{\mathcal{P}})\in|{\bf LFCInt}|$, where $i_{\mathcal{P}}:L^X\rightarrow L^X$ is a \v{C}ech $L$-valued fuzzy interior operator such that for each $x\in X,f\in L^X,$
\[i^{{\mathcal{P}}}(f)(x)=\bigwedge\limits_{x'\in X}\neg A_{\xi(x)}(x')\#f(x')=F_{\mathcal{P}}^{\downarrow}[f](\xi(x)),\] where $\xi(x)\in J$ is the unique index such that $x\in core(A_{\xi(x)})$.
\end{itemize}
\begin{pro}\label{23}
If $(\phi,\mu):\pi_1\rightarrow \pi_2$ is a ${\bf LSpaceFP}$-morphism, then  $(\varrho,\varrho'):\Sigma_{\pi_1}\rightarrow {\Sigma_{\pi_2}}$ is ${\bf LFPrTop}$-morphism, where $ \varrho:X_1\times X_2\rightarrow L$ is an $L$-valued fuzzy relation such that $\varrho(x_1,x_2)=\phi(x_1,x_2),\,\forall\,x_1\in X_1,x_2\in X_2$ and $\varrho':X_2\times X_1\rightarrow L$ is also an $L$-valued fuzzy relation such that $\varrho'(x_2,x_1)=\mu(\xi'(x_2),\xi(x_1)),\,\forall\,x_2\in core(A_{\xi'(x_2)}),x_1\in core(A_{\xi(x_1)}),\xi'(x_2)\in J_2,\xi(x_1)\in J_1$.
\end{pro}
\textbf{Proof:} Let $(\phi,\mu):\pi_1\rightarrow \pi_2$ be an ${\bf LSpaceFP}$-morphism. Then  for all $x_1\in X,f\in L^{X_2}$ 
\begin{eqnarray*}
p^{{\mathcal{P}_1}}_{x_1}(\overleftarrow{\varrho}(f))&=&F_{\mathcal{P}_1}^{\downarrow}[\overleftarrow{\varrho}(f)](\xi(x_1))\\
&=&\bigwedge\limits_{x_1'\in X_1}\neg A_{\xi(x_1)}(x_1')\#\overleftarrow{\varrho}(f)(x_1')\\
&=&\bigwedge\limits_{x_1'\in X_1}\neg A_{\xi(x_1)}(x_1')\#\bigwedge\limits_{x'_2\in X_2}\neg\varrho(x_1',x'_2)\#f(x_2')\\
&=&\bigwedge\limits_{x_1'\in X_1}\bigwedge\limits_{x'_2\in X_2}\neg (A_{\xi(x_1)}(x_1')\ast\varrho(x_1',x'_2))\#f(x_2')\\
&\geq&\bigwedge\limits_{x'_2\in X_2}\neg (\mu(\xi'(x_2),\xi(x_1))\#A'_{\xi'(x_2)}(x'_2))\#f(x_2')\\
&=&\bigwedge\limits_{x'_2\in X_2}(\neg \mu(\xi'(x_2),\xi(x_1))\ast \neg A'_{\xi'(x_2)}(x'_2))\#f(x_2')\\
&\geq&\bigwedge\limits_{x'_2\in X_2}\neg \mu(\xi'(x_2),\xi(x_1))\ast (\neg A'_{\xi'(x_2)}(x'_2)\#f(x_2'))\\
&\geq&\neg \mu(\xi'(x_2),\xi(x_1))\ast \bigwedge\limits_{x'_2\in X_2} (\neg A'_{\xi'(x_2)}(x'_2)\#f(x_2'))\\
&=&\neg \varrho'(x_2,x_1)\ast F_{\mathcal{P}_2}^{\downarrow}[f](\xi'(x_2))\\
&=&\neg \varrho'(x_2,x_1)\ast p^{{\mathcal{P}_2}}_{x_2}(f).
\end{eqnarray*}
Thus $\neg \varrho'(x_2,x_1)\ast p^{{\mathcal{P}_2}}_{x_2}(f)\leq p^{{\mathcal{P}}}_{x_1}(\overleftarrow{\varrho}(f))$. Therefore
$(\varrho,\varrho'):\Sigma_{\pi_1}\rightarrow \Sigma_{\pi_2}$ is an ${\bf LFPrTop}$-morphism.
\begin{pro}\label{F_7}
Let $F_7:{\bf LSpaceFP}\rightarrow {\bf LFPrTop}$ be a function such that for every $\pi_1\in |{\bf LSpaceFP}|, F_7(\pi_1)={\Sigma_{\pi_1}}$ and for every ${\bf LSpaceFP}$-morphism $(\phi,\mu):\pi\rightarrow\pi_1,F_7(\phi,\mu):{\Sigma_{\pi_1}}\rightarrow {\Sigma_{\pi_2}}$ be a function such that $F_7(\phi,\mu)=(\varrho,\varrho')$. Then $F_7$ is a functor.
\end{pro}
\textbf{Proof:} (i) Let $(\phi_1,\mu_1):\pi_1\rightarrow\pi_2,(\phi_2,\mu_2):\pi_1\rightarrow\pi_2$ be  ${\bf LSpaceFP}$-morphisms. Then 
\begin{eqnarray*}
F_7((\phi_2,\mu_2)\bullet(\phi_1,\mu_1))&=&(\varrho_2,\varrho'_2)\bullet(\varrho_1,\varrho'_1)\\
&=&F_7(\phi_2,\mu_2)\bullet F_7(\phi_2,\mu_2).
\end{eqnarray*}
Therefore $F_7((\phi_2,\mu_2)\bullet(\phi_1,\mu_1))=F_7(\phi_2,\mu_2)\bullet F_7(\phi_1,\mu_1)$.\\\\
(ii) Let $\pi_1=(X_1,J_1,\mathcal{P}_1)\in|{\bf LSpaceFP}|$. Then $F_7(id_{\pi_1})=F_7(1_X,0_J)=(1_X,0_X)$. Thus $F_7(id_{\pi_1})=id_{F_7({\pi_1})}$. Hence $F_7$ is a functor.
\begin{pro}\label{F_8}
Let $F_8:{\bf LSFtrans}^\downarrow\rightarrow {\bf LFPrTop}$ be a function such that for every $\mathcal{H}_1\in |{\bf LSFtrans}^\downarrow|, F_8(\mathcal{H}_1)={\Sigma_{\mathcal{H}_1}}$ and for every ${\bf LFtrans}^\downarrow$-morphism $(\psi,\nu):\mathcal{H}_1\rightarrow\mathcal{H}_2,F_8(\psi,\nu):{\Sigma_{\mathcal{H}_1}}\rightarrow {\Sigma_{\mathcal{H}_2}}$ be a function such that $F_8(\psi,\nu)=(\varrho,\varrho')$. Then $F_8$ is a functor.
\end{pro}
\textbf{Proof:} {Propositions \ref{16} and \ref{F_7} lead to this proof.}
\begin{pro}\label{00}
If $(\phi,\mu):\pi_1\rightarrow \pi_2$ is a ${\bf LSpaceFP}$-morphism, then  $(\kappa,\kappa'):I_{\pi}\rightarrow {I_{\pi_1}}$ is ${\bf LFCInt}$-morphism, where $ \kappa:X_1\times X_2\rightarrow L$ is an $L$-valued fuzzy relation such that $\kappa(x_1,x_2)=\phi(x_1,x_2),\,\forall\,x_1\in X_1,x_2\in X_2$ and $\kappa':X_2\times X_1\rightarrow L$ is also an $L$-valued fuzzy relation such that $\kappa'(x_2,x_1)=\mu(\xi'(x_2),\xi(x_1)),\,\forall\,x_2\in core(A_{\xi'(x_2)}),x_1\in core(A_{\xi(x_1)}),\xi'(x_2)\in J_2,\xi(x_1)\in J_1$.
\end{pro}
\textbf{Proof:} Let $(\phi,\mu):\pi\rightarrow \pi_1$ be an ${\bf LSpaceFP}$-morphism. Then  for all $x_1\in X,f\in L^{X_2}$ 
\begin{eqnarray*}
i^{{\mathcal{P}_1}}(\overleftarrow{\kappa}(f))(x_1)&=&F_{\mathcal{P}_1}^{\downarrow}[\overleftarrow{\kappa}(f)](\xi(x_1))\\
&=&\bigwedge\limits_{x_1'\in X}\neg A_{\xi(x_1)}(x_1')\#\overleftarrow{\kappa}(f)(x_1')\\
&=&\bigwedge\limits_{x_1'\in X}\neg A_{\xi(x_1)}(x_1')\#\bigwedge\limits_{x'_2\in X_2}\neg\kappa(x_1',x'_2)\#f(x_2')\\
&=&\bigwedge\limits_{x_1'\in X}\bigwedge\limits_{x'_2\in X_2}\neg (A_{\xi(x_1)}(x_1')\ast\kappa(x_1',x'_2))\#f(x_2')\\
&\geq&\bigwedge\limits_{x'_2\in X_2}\neg (\mu(\xi'(x_2),\xi(x_1))\#A'_{\xi'(x_2)}(x'_2))\#f(x_2')\\
&=&\bigwedge\limits_{x'_2\in X_2}(\neg \mu(\xi'(x_2),\xi(x_1))\ast \neg A'_{\xi'(x_2)}(x'_2))\#f(x_2')\\
&\geq&\bigwedge\limits_{x'_2\in X_2}\neg \mu(\xi'(x_2),\xi(x_1))\ast (\neg A'_{\xi'(x_2)}(x'_2)\#f(x_2'))\\
&\geq&\neg \mu(\xi'(x_2),\xi(x_1))\ast \bigwedge\limits_{x'_2\in X_2} (\neg A'_{\xi'(x_2)}(x'_2)\#f(x_2'))\\
&=&\neg \kappa'(x_2,x_1)\ast F^\downarrow_{{\mathcal{P}_2}}[f]({\xi'(x_2)})\\
&=&\neg \kappa'(x_2,x_1)\ast i^{{\mathcal{P}_2}}(f)({x_2}).
\end{eqnarray*}
Thus $\neg \kappa'(x_2,x_1)\ast i^{{\mathcal{P}_2}}(f)({x_2})\leq i^{{\mathcal{P}_1}}(\overleftarrow{\kappa}(f))({x_1})$. Therefore
$(\kappa,\kappa'):I_{\pi_1}\rightarrow I_{\pi_2}$ is an ${\bf LFCInt}$-morphism.
\begin{pro}
Let $F_9:{\bf LSpaceFP}\rightarrow {\bf LFCInt}$ be a function such that for every $\pi_1\in |{\bf LSpaceFP}|, F_9(\pi_1)={I_{\pi_1}}$ and for every ${\bf LSpaceFP}$-morphism $(\phi,\mu):\pi-1\rightarrow\pi_2,F_9(\phi,\mu):I_{\pi_1}\rightarrow {I_{\pi_2}}$ be a function such that $F_9(\phi,\mu)=(\kappa,\kappa')$. Then $F_9$ is a functor.
\end{pro}
\textbf{Proof:} Similar to that of Proposition \ref{F_7}.
\begin{pro}
Let $F_{10}:{\bf LSFtrans}^\downarrow\rightarrow {\bf LFCInt}$ be a function such that for every $\mathcal{H}_1\in |{\bf LSFtrans}^\downarrow|, F_{10}(\mathcal{H}_1)=I_{\mathcal{H}_1}$ and for every ${\bf LFtrans}^\downarrow$-morphism $(\psi,\nu):\mathcal{H}_1\rightarrow\mathcal{H}_2,F_{10}(\psi,\nu):I_{\mathcal{H}_1}\rightarrow {I_{\mathcal{H}_2}}$ be a function such that $F_{10}(\psi,\nu)=(\kappa,\kappa')$. Then $F_{10}$ is a functor.
\end{pro}
\textbf{Proof:} Similar to that of Proposition \ref{F_8}.\\\\
Next, we have the following.
\begin{pro}\label{com} Let $F_3:{\bf LSpaceFP}\rightarrow {\bf LFtrans}^{\downarrow},F'_3:{\bf LFtrans}^{\downarrow}\rightarrow {\bf LSpaceFP},F_6:{\bf LFPrTop}\rightarrow{\bf LFCInt},F_7:{\bf LSpaceFP}\rightarrow {\bf LFPrTop},F_8:{\bf LFtrans}^{\downarrow}\rightarrow {\bf LFPrTop},F_9:{\bf LSpaceFP}\rightarrow {\bf LFCInt}$ and $F_{10}:{\bf LFtrans}^{\downarrow}\rightarrow{\bf LFCInt}$ be functors. Then the diagram in Figure \ref{fig:4} is commutable.
\begin{figure}
\begin{center}
\begin{tikzcd}[row sep=10ex, column sep=10ex] 
{\bf LSpaceFP}{{\arrow{r}{F_7}}}\arrow[swap]{d}{F_3}& {\bf LFPrTop} \arrow{d}{F_6}\arrow[r, "F_8", leftarrow] & {\bf LFtrans}^{\downarrow}\arrow{d}{F_3'}\\
{\bf LFtrans}^{\downarrow}\arrow{r}{F_{10}}& {\bf LFCInt}\arrow[r, "{F_9}", leftarrow] &{\bf LSpaceFP}
\end{tikzcd}
\end{center}
\caption{Diagram for Proposition \ref{com}.}
\label{fig:4}
  \end{figure}
\end{pro}
\textbf{Proof:} Let $\pi_1\in |{\bf LSpaceFP}|$ and $(\phi,\mu):\pi_1\rightarrow \pi_2$ be an ${\bf LSpaceFP}$-morphism. Then 
\begin{eqnarray*}
(F_{10}\circ F_3)(\pi)&=&F_{10}(F_3(\pi_1))\\
 &=&F_{10}(\mathcal{H}_{\pi_1})\\
 &=&I_{\mathcal{H}_{\pi_1}}\\
 &=&I_{\Sigma_{\pi_1}}\\
(F_6\circ F_7)(\pi_1)&=&F_6(F_7(\pi_1))\\
&=&F_6(\Sigma_{\pi_1})\\
&=&I_{\Sigma_{\pi_1}}.
\end{eqnarray*}
Thus $(F_{10}\circ F_3)(\pi_1)=(F_6\circ F_7)(\pi_1)$ and 
\begin{eqnarray*}
 (F_{10}\circ F_3)(\phi,\mu)&=&F_{10}(F_3(\phi,\mu))\\
 &=&F_{10}(\psi,\nu)\\
 &=&(\kappa,\kappa')\\
 (F_6\circ F_7)(\phi,\mu)&=&F_6(F_7(\phi,\mu))\\
 &=&F_6(\varrho,\varrho')\\
 &=&(\kappa,\kappa').
\end{eqnarray*}
Thus $(F_{10}\circ F_3)(\phi,\mu)=(F_6\circ F_7)(\phi,\mu)$. Therefore $F_{10}\circ F_3=F_6\circ F_7$. Similarly, for every $\mathcal{H}_1\in |{\bf LFtrans}^\downarrow|$ and for every ${\bf LFtrans}^\downarrow$-morphism $(\psi,\nu):\mathcal{H}_1\rightarrow \mathcal{H}_2$, we can show that $F_{6}\circ F_8=F_9\circ {F_3'}$. Thus diagram in Figure \ref{fig:4} is commutable.\\\\
Finally, we demonstrate the existence of two pairs of functors $(F_3,F'_3)$ and $(F_6,F'_6)$ that have adjoint property. Now, we begin with following.
\begin{pro}\label{adj}
Let $F_3:{\bf LSpaceFP}\rightarrow {\bf LFtrans}^\downarrow$ and $F'_3:{\bf LFtrans}^\downarrow\rightarrow {\bf LSpaceFP}$ be functors. Then $F_3$ is left adjoint to $F'_3$ and $F'_3$ is a right adjoint to $F_3$.
\end{pro}
\textbf{Proof:} To show this result, we have to show that there exists a natural transformation $\Psi:id_{\bf LSpaceFP}\rightarrow F'_3\circ F_3$ such that for every $\pi\in |{\bf LSpaceFP}|$ and ${\bf LSpaceFP}$-morphism $(\phi,\mu):\pi\rightarrow F_3'(\mathcal{H})$, there exists a unique ${\bf LFtrans}^\downarrow$-morphism $(\psi,\nu):F_3(\pi)\rightarrow \mathcal{H}$ such that the diagram in Figure \ref{fig:6} commutes. To do this, let $\pi=(X,J,\mathcal{P})\in |{\bf LSpaceFP}|,\mathcal{H}=(X,J,v,H)\in|{\bf LFtrans}^\downarrow|$ and $\Psi_\pi:(X,J,\mathcal{P})\rightarrow (F'_3\circ F_3)(X,J,\mathcal{P})$ be a function such that $\Psi_\pi=(1_X,0_J)$. It is easy to check that $\Psi$ is a natural transformation. Further, let $(\phi,\mu):(X,J,\mathcal{P})\rightarrow F'_3(X,J,v,H)$ be a ${\bf LSpaceFP}$-morphism. Now, we define an ${\bf LFtrans}^\downarrow$-morphism $(\psi,\nu):F_3(X,J,\mathcal{P})\linebreak\rightarrow (X,J,v,H)$ such that $(\psi,\nu)=(\phi,\mu)$. The diagram in Figure \ref{fig:7} commutes, i.e.,
\begin{eqnarray*}
F_3'(\psi,\nu)\bullet (1_X,0_J)&=&(\psi,\nu)\bullet (1_X,0_J)\\
&=&(\psi\odot 1_X,\nu\oplus 0_J)\\
&=&(\psi,\nu)\\
&=&(\phi,\mu).
\end{eqnarray*}
Thus $F_3'(\psi,\nu)\bullet (1_X,0_J)=(\phi,\mu)$. The uniqueness of $(\psi,\nu)$ is obvious. Thus $F_3$ is a left adjoint of $F_3'$ and $F_3'$ is a right adjoint of $F_3$.
\begin{figure}
\[
 \begin{tikzcd}[row sep=12ex, column sep=12ex] 
   \pi\arrow{r}{\Psi_\pi} \arrow[swap]{dr}{(\phi,\mu)} & F'_3(F_3(\pi))\arrow{d}{F'_3(\psi,\nu)} &F_3(\pi)\arrow{d}{(\psi,\nu)}\\
     & F'_3(\mathcal{H})&\mathcal{H}
  \end{tikzcd}
  \]
 \caption{Diagram for Proposition \ref{adj}.}
\label{fig:6}
  \end{figure}
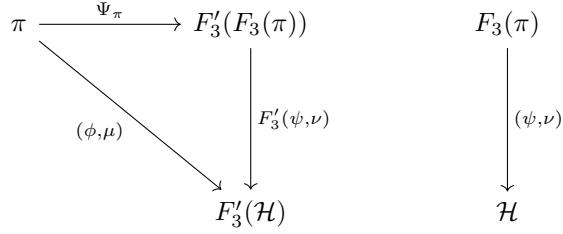
  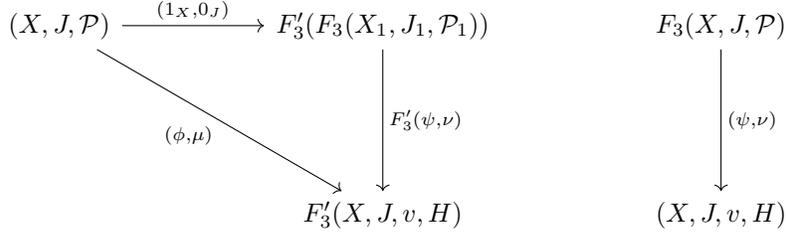
\begin{figure}
  \[
  \begin{tikzcd} [row sep=12ex, column sep=12ex] 
  \
    (X,J,\mathcal{P}) \arrow{r}{(1_X,0_J)} \arrow[swap]{dr}{(\phi,\mu)} & F'_3(F_3(X_1,J_1,\mathcal{P}_1))\arrow{d}{F'_3(\psi,\nu)} &F_3(X,J,\mathcal{P})\arrow{d}{(\psi,\nu)}\\
     &  F'_3(X,J,v,H)&(X,J,v,H)
    \end{tikzcd}
\]
\caption{Diagram is equivalent to the diagram in Figure \ref{fig:6}.}
\label{fig:7}
  \end{figure}
\begin{pro}\label{adj1}
Let $F_6:{\bf LFPrTop}\rightarrow {\bf LFCInt}$ and $F'_6:{\bf LFCInt}\rightarrow {\bf LFPrTop}$ be functors. Then $F_6$ is a left adjoint to $F'_6$ and $F'_6$ is a right adjoint to $F_6$.
\end{pro}
\textbf{Proof:} To show this result, we have to show that there exists a natural transformation $\Theta:id_{\bf LFPrTop}\rightarrow F'_6\circ F_6$ such that for every $\Sigma\in |{\bf LFPrTop}|$ and ${\bf LFPrTop}$-morphism $(\varrho,\varrho'):\Sigma\rightarrow F_6'(I)$, there is a unique ${\bf LFCInt}$-morphism $(\kappa,\kappa'):F_6(\Sigma)\rightarrow I$ such that the diagram in Figure \ref{fig:8} commutes. To do this, let $\Sigma=(X,P)\in |{\bf LFPrTop}|,I=(X,i)\in|{\bf LFCInt}|$ and $\Theta_\Sigma:(X,{P})\rightarrow (F'_6\circ F_6)(X,{P})$ be a function such that $\Theta_\Sigma=(1_X,0_X)$. It is easy to check that $\kappa$ is a natural transformation. Also, let $(\varrho,\varrho'):(X,{P})\rightarrow F'_6(X,i)$ be a ${\bf LFPrTop}$-morphism. Now, we define an ${\bf LFCInt}$-morphism $(\kappa,\kappa'):F_6(X,P)\rightarrow (X,i)$ such that $(\kappa,\kappa')=(\varrho,\varrho')$. The diagram in Figure \ref{fig:9} commutes, i.e.,
\begin{eqnarray*}
F_6'(\kappa,\kappa')\bullet (1_X,0_X)&=&(\kappa,\kappa')\bullet (1_X,0_X)\\
&=&(\kappa\odot 1_X,\kappa'{\oplus} 0_X)\\
&=&(\kappa,\kappa')\\
&=&(\varrho,\varrho').
\end{eqnarray*}
Thus $F_6'(\kappa,\kappa')\bullet (1_X,0_X)=(\varrho,\varrho')$. The uniqueness of $(\kappa,\kappa')$ is obvious. Thus $F_6$ is left adjoint of $F_6'$ and $F_6'$ is a right adjoint of $F_6$.
\begin{figure}
\[
 \begin{tikzcd}[row sep=12ex, column sep=12ex] 
   \Sigma\arrow{r}{\Theta_\Sigma} \arrow[swap]{dr}{(\varrho,\varrho')} & F'_6(F_6(\Sigma))\arrow{d}{F'_6(\kappa,\kappa')} &F_6(\Sigma)\arrow{d}{(\kappa,\kappa')}\\
     & F'_6(I)&I
  \end{tikzcd}
  \]
 \caption{Diagram for Proposition \ref{adj1}.}
\label{fig:8}
  \end{figure}
  \begin{figure}
    \[
  \begin{tikzcd}[row sep=12ex, column sep=10ex]  
  (X,P)\arrow{r}{(1_X,0_X)} \arrow[swap]{dr}{(\varrho,\varrho')} & F'_6(F_6(X,P))\arrow{d}{F'_6(\kappa,\kappa')} &F_6(X,P)\arrow{d}{(\kappa,\kappa')}\\
     & F'_6(X,i)&(X,i)
    \end{tikzcd}
\]
\caption{Diagram is equivalent to the diagram in Figure \ref{fig:8}.}
\label{fig:9}
  \end{figure}
  \section{Conclusion}
In view of the work carried out in \cite{shai} towards the importance of ${\bf Qua}$ category in the unification of categories associated with fuzzy automata, herein, we have tried to unify the categories associated with $F$-transforms and fuzzy pretopological spaces. Specifically, we introduced the categories ${\bf LSpaceFP}$, ${\bf LFtrans}^\downarrow$, ${\bf LFPrTop}$ and ${\bf LFCInt}$ and their morphisms are as pairs of $L$-valued fuzzy relations (respectively, the categories ${\bf SpaceFP}$, ${\bf Ftrans}^\downarrow$, ${\bf FPrTop}$ and ${\bf CInt}$ and their morphisms are as functions). Interestingly, we have underlaid the connections of the categories ${\bf LSpaceFP}$, ${\bf LFtrans}^\downarrow$, ${\bf LFPrTop}$ and ${\bf LFCInt}$ with the category ${\bf Qua}$. Further, we show a functorial relationship among the categories ${\bf LSpaceFP}$, ${\bf LFtrans}^\downarrow$, ${\bf LFPrTop}$ and ${\bf LFCInt}$. Interestingly,  there exist two pairs of functors $(F_3,F'_3)$ and $(F_6,F'_6)$ that have adjoint property. In the future, it will be interesting to see the unification of categories associated with the general fuzzy relations in the framework of the ${\bf Qua}$ category.
\section*{Compliance with ethical standards}
{\bf Ethical approval} This paper does not deal with any ethical problems.\\\\
{\bf Funding details} This work did not receive any grant from funding agencies
in the public, commercial, or not-for-profit sectors.\\\\
{\bf Conflict of interest} The authors declared that they have no conflicts of interest to this work.\\\\
{\bf Informed consent} We declare that all authors have informed Consent.\\\\
{\bf Author contributions} Abha Tripathi was involved in conceptualization, methodology, and writing-original draft preparation; S.P. Tiwari was involved in supervision and writing-review and editing. All authors attest that they meet the criteria for authorship.


\begin{thebibliography}{99}
\bibitem{ad} J. Ad\'{a}mek, H. Herlich, G.E. Strecker, Abstract and concrete categories, {\em Wiley, New York}, (1990).
\bibitem{ar} M.A. Arbib, E.G. Manes, Arrows, structures, and functors: the categorical imperative, {\em Academic Press: New York}, (1975).
\bibitem{bar1} M. Barr, Fuzzy models of linear logic, {\em Mathematical Structures in Computer Science}, {\bf 6} (1996) 301-312.
\bibitem{bar} M. Barr, C. Wells, Category theory for computing science, {\em New York: Prentice Hall}, (1990).
\bibitem {belo} R. B\v{e}lohl\'{a}vek, Fuzzy relational systems, {\em  Kluwer Academic Publishers, Plenum Publishers}, New York, Boston, Dordrecht, London, Moscow, (2002).
\bibitem{blas} A. Blass, A category arising in linear logic, complexity theory and set theory, {\em Advances in Linear Logic}, (1995) 61-81.
\bibitem{pai} V. De Paiva, The Dialectica categories, {\em In Proceedings of the Conference on Categories in Computer Science and Logic}, (1989) 4-62.
\bibitem{mar} F. Di Martino, V. Loia, I. Perfilieva, S. Sessa, An image coding/decoding method based on direct and inverse fuzzy transforms, {\em International Journal of Approximate Reasoning}, {\bf 48} (2008) 110-131.
{ \bibitem{mar1} F. Di Martino, V. Loia, S. Sessa, A segmentation method for images compressed by fuzzy transforms, {\em Fuzzy Sets and Systems}, {\bf 161} (2010) 56-74.}
\bibitem{ma} F. Di Martino, V. Loia, S. Sessa, Fuzzy transforms method in prediction data analysis, {\em Fuzzy Sets and Systems}, {\bf 180} (2011) 146-163.
\bibitem{eil} S. Eilenberg, S.M. Lane, General theory of natural equivalences, {\em Transactions of the American Mathematical Society}, {\bf 58} (1945) 231-294.
\bibitem{fre} P.J. Freyd, Abelian categories, {\em New York: Harper and Row}, (1964).
\bibitem{gog} J.A. Goguen, $L$-fuzzy sets, \emph{Journal of Mathematical Analysis and Applications}, {\bf18} (1967) 147-174.
\bibitem{gro} A. Grothendieck, Sur quelques points d'alg\'{e}bre homologique, {\em Tohoku
Mathematical Journal}, {\bf 9} (1957) 119-183.
\bibitem{holc} M. Hol\v{c}apek, L. Nguyen, Trend-cycle estimation using fuzzy transform of higher degree, {\em Iranian Journal of Fuzzy Systems}, {\bf15} (2018) 23-54.
{\bibitem{hut} P. Hurt\'{\i}k, S. Tomasiello, A review on the application of fuzzy transform in data and image compression, {\em Soft Computing}, {\bf 23} (2019) 12641-12653.}
\bibitem{kh} A. Khastan, I. Perfilieva, Z. Alizani, A new fuzzy approximation method to Cauchy problems by fuzzy transform, {\em Fuzzy Sets and Systems,} {\bf 288} (2016) 75-95.
{\bibitem{kh1} A. Khastan, {A new representation for inverse fuzzy transform and its application}, {\em Soft Computing}, {\bf 21} (2017), 3503-3512.}
\bibitem{kl} E.P. Klement, R. Mesiar, E. Pap, Triangular norms, trends in logic, {\em Kluwer Academic Publishers}, Dordrecht, {\bf8} (2000).
\bibitem{kli} G.J. Klir, B. Yuan, Fuzzy logic: theory and applications, {\em Prentice Hall}, Englewood Cliffs, NJ, (1995).
\bibitem{laf} Y. Lafont and T. Streicher, Game semantics for linear logic, {\em In Proceedings of $6^{th}$ Annual IEEE Symposium on Logic in Computer Science}, (1991) 43-49.
\bibitem{lan} S.M. Lane, Categories for the working mathematician, {\em Springer Science and Business Media}, (2013).
\bibitem{law} F.W. Lawvere, Functorial semantics of algebraic theories, {\em In Proceedings of the National Academy of Sciences of the United States of America}, {\bf 50} (1963) 1-869.
\bibitem{law1} F.W. Lawvere, The category of categories as a foundation for mathematics, In: {\em Proceedings of the Conference on Categorical Algebra}, (1966) 1-20.
\bibitem{li} M. Liu, D. Chen, C. Wu, H. Li, Approximation theorem of the fuzzy transform in fuzzy reasoning and its application to the scheduling problem, {\em Computers and Mathematics with Applications}, {\bf 51} (2006) 515-526.
\bibitem{mockor1} J. Mo\v{c}ko\v{r}, Fuzzy type relations and transformation operators defined by monads, {International Journal of Computational Intelligence Systems}, {\bf 13} (2020) 1530-1538.
\bibitem{mock} J. Mo\v{c}ko\v{r}, Spaces with fuzzy partitions and fuzzy transform, {\em Soft Computing}
{\bf 13} (2017), 3479-3492.
\bibitem{jiri} J. Mo\v{c}ko\v{r}, Axiomatic of lattice-valued $F$-transform, {\em Fuzzy Sets and Systems}, {\bf342}
(2018) 53-66.
\bibitem{jir} J. Mo\v{c}ko\v{r}, $F$-transforms and semimodule homomorphisms, {\em Soft Computing}, {\bf 23} (2019) 7603-7619.
\bibitem{jiri1} J. Mo\v{c}ko\v{r}, Relational, closure and partition powerset theories, {\em Fuzzy Sets and Systems}, {\bf 420} (2021) 100-122.
\bibitem{mo} J. Mo\v{c}ko\v{r}, M. Hol\v{c}apek, Fuzzy objects in spaces with fuzzy partitions, {\em Soft Computing}, {\bf 21} (2016) 7268-7284.
\bibitem{mocko} J. Mo\v{c}ko\v{r}, P. Hurt\'{\i}k, Lattice-valued $F$-transforms and similarity relations, {\em Fuzzy Sets and Systems}, {\bf 342} (2018) 67-89.
\bibitem{mockor} J. M\v{o}ck\v{o}r, I. Perfilieva, Functors among categories of $L$-fuzzy partitions, $L$-fuzzy pretopological spaces and $L$-fuzzy closure spaces {\em In IFSA World Congress and NAFIPS Annual Conference: Fuzzy Techniques: Theory and Applications, IFSA/NAFIPS Louisiana, USA}, Springer: Cham, Switzerland, (2019) 382-393.
\bibitem{vil} V. Nov\'{a}k, I. Perfilieva, M. Hol\v{c}apek, V. Kreinovich, {Filtering out high frequencies in time series using F-transform}, {\em Information Sciences}, {\bf 274} (2014) 192-209.
\bibitem{per} I. Perfilieva, $F$-transforms: theory and its applications, \emph{Fuzzy Sets and Systems}, {\bf157} (2006) 993-1023.
\bibitem{irin1} I. Perfilieva, Fuzzy transforms: a challenge to conventional  transforms, {\em Advances in Image and Electron Physics}, {\bf 147} (2007) 137-196.
\bibitem{no} I. Perfilieva, V. Nov\'{a}k, A. {Dvo\v{r}\'{a}k}, {Fuzzy transforms in the analysis of data}, {\em International Journal of Approximate Reasoning}, {\bf 48} (2008) 36-46.
\bibitem{anan} I. Perfilieva, A.P. Singh, S.P. Tiwari, On the relationship among $F$-transform, fuzzy rough sets and fuzzy topology, \emph{Soft Computing}, {\bf21} (2017) 3513-3523.
\bibitem{spt1} I. Perfilieva, S.P. Tiwari, A.P. Singh, Lattice-valued $F$-transforms as interior operators of $L$-fuzzy pretopological spaces, {\em Communications in Computer and Information Science}, {\bf 854} (2018) 163-174.
\bibitem{Ir} I. Perfilieva, R. Valasek, {Fuzzy transforms in removing noise}, {\em Advances in Soft Computing}, {\bf 2} (2005) 221-230.
\bibitem{pi} B.C. Pierce, Basic Category Theory for Computer Scientists, {\em The MIT Press, Cambridge}, (1991).
\bibitem{qi} J. Qiao, B.Q. Hu, A short note on $L$-fuzzy approximation spaces and $L$-fuzzy pretopological spaces, {\em Fuzzy Sets and Systems}, {\bf 312} (2017) 126-134. 
\bibitem{qia} J. Qiao, B.Q. Hu, On $(\otimes,\&)$-fuzzy rough sets based on residuated and co-residuated lattice, {\em Fuzzy Sets and Systems}, {\bf 336} (2017) 54-86. 
\bibitem{roh} S. B. Roh, S. K. Oh, J. H. Yoon, K. Seo, {Design of face recognition system based on fuzzy transform and radial basis fnction neural networks}, {\em Soft Computing}, {\bf 23} (2019) 4969-4985.
\bibitem{shai} S. Singh, S.P. Tiwari, On unification of categories of fuzzy automata as Qua category, {\em Soft Computing}, {\bf 26} (2022) 1509-1529.
\bibitem{sinha} P. Sinha, Algebraic nondeterministic and transition systems, Ph.D. Thesis, IIT Delhi, (2005).
\bibitem {st} L. Stefanini, $F$-transform with parametric generalized fuzzy partitions, {\em Fuzzy Sets and Systems}, {\bf 180} (2011) 98-120.
\bibitem {step} M. {\v{S}t\v{e}pni\v{c}ka}, R. {Val\'{a}\v{s}ek}, Fuzzy transforms and their application to wave equation, {\em Journal of Electrical Engineering}, {\bf 55} (2004) 7-10.
\bibitem {ste} M. {\v{S}t\v{e}pni\v{c}ka}, O. {Polakovi\v{c}}, {A neural network approach to the fuzzy transform}, {Fuzzy Sets and Systems}, {\bf 160} (2009) 1037-1047.
\bibitem{tri} A. Tripathi, S.P. Tiwari, I. Perfilieva, { $F$-transforms determined by implicators}, {\em Iranian Journal of Fuzzy Systems}, {\bf 18} (2021) 19-36.
\bibitem{abhaog} A. Tripathi, S.P. Tripathi,  Kavikumar Jacob, S. Mahto, $F$-transforms determined by overlap and grouping maps over a complete lattice, {\em Soft Computing}, {28} 2024 10781-10800.
\bibitem{abha} A. Tripathi, S.P. Tripathi,  Kavikumar Jacob, D. Nagarajan, A fuzzy function granular 
$F$-transform and inverse-transform with application, {\em Decision Analytics Journal}, {7} 2023 100241.
\bibitem{abhac} {A. Tripathi}, S.P. Tiwari, A.P. Singh, On $L^M$-valued $F$-transforms, $L^M$-valued fuzzy rough sets, \textit{In Proceedings of the Conference of the International Fuzzy Systems Association and the European Society for Fuzzy Logic and Technology (EUSFLAT 2019)}, \textbf{1} (2019) 220-226. 
\bibitem{trip} A. Tripathi, S.P. Tiwari, A.P. Singh, {On $L^M$-valued $F$-transforms, $L^M$-valued fuzzy rough sets and $L^M$-valued fuzzy transformation systems}, {\em New Mathematics and Natural Computation}, {\bf 17} (2021) 339-359.
\bibitem{to} L. Troiano, P. Kriplani, {Supporting trading strategies by inverse fuzzy transform}, {\em Fuzzy Sets and Systems}, {\bf 180} (2011) 121-145.
\bibitem{zh} D. Zhang, Fuzzy pretopological spaces, an extensional topological extension of FTS, {\em Chinese Annals of Mathematics}, {\bf 3} (1999) 309-316. 
\end{thebibliography}
\end{document}